\documentclass[mathpazo]{cicp}

\usepackage{amsmath}
\usepackage{amsthm}
\usepackage{newtxmath}
\usepackage[T1]{fontenc}
\usepackage[utf8]{inputenc}
\usepackage{units}
\usepackage{graphicx}
\usepackage{float}
\usepackage{algorithmic}
\usepackage{framed,multirow}
\usepackage{amssymb}
\usepackage{latexsym}
\usepackage{url}
\usepackage{graphics}
\usepackage{graphicx}
\usepackage{enumerate}
\usepackage{subfigure}
\usepackage{float}
\usepackage{colortbl}
\usepackage{color,xcolor}
\definecolor{light-gray}{gray}{0.95} 
\usepackage{booktabs}% ?    \toprule  \midrule  \bottomrule
\usepackage{graphbox}
\usepackage{mathrsfs}
\usepackage{doi}
\usepackage{hyperref}
\usepackage{caption}
\usepackage[ruled,linesnumbered]{algorithm2e}
\definecolor{newcolor}{rgb}{.8,.349,.1}
\usepackage{changepage}
\usepackage{amsfonts}% blackboard math symbols
\usepackage{nicefrac}% compact symbols for 1/2, etc.
\begin{document}

\title{Solving multiscale elliptic problems by sparse radial basis function neural networks}

\author[Zhiwen]{Zhiwen Wang\affil{1}, Minxin Chen\affil{1}, and Jingrun Chen\affil{2}\comma\corrauth}
\address{
	\affilnum{1}\ School of Mathematical Sciences, Soochow University, Suzhou and 215006, China \\
	\affilnum{2}\ School of Mathematical Sciences and Suzhou Institute for Advanced Research, University of Science and Technology of China, Hefei and 230026, China}
\emails{{\tt 20224007006@stu.suda.edu.cn} (Z. Wang), {\tt chenminxin@suda.edu.cn} (M. Chen), {\tt jingrunchen@ustc.edu.cn} (J. Chen).}

%\title{Solving multiscale elliptic problems by sparse radial basis function neural networks}%
%
%\author[1]{Zhiwen \snm{Wang}}
%\email{20224007006@stu.suda.edu.cn}
%\author[1]{Minxin \snm{Chen}\corref{cor1}}
%\ead{chenminxin@suda.edu.cn}
%\author[2]{Jingrun \snm{Chen}\corref{cor2}}
%\ead{jingrunchen@ustc.edu.cn}
%
%\cortext[cor1]{Corresponding author: School of Mathematical Sciences, Soochow University, Suzhou, 215006, China.}
%\cortext[cor2]{Corresponding author: School of Mathematical Sciences and Suzhou Institute for Advanced Research, University of Science and Technology of China, Hefei and 230026, China.}
%
%\address[1]{School of Mathematical Sciences, Soochow University, Suzhou and 215006, China}
%\address[2]{School of Mathematical Sciences and Suzhou Institute for Advanced Research, University of Science and Technology of China, Hefei and 230026, China}

\begin{abstract}
%%%
Machine learning has been successfully applied to various fields of scientific computing in recent years. For problems with multiscale features such as flows in porous media and mechanical properties of composite materials, however, it remains difficult for machine learning methods to resolve such multiscale features, especially when the small scale information is present. In this work, we propose a sparse radial basis function neural network method to solve elliptic partial differential equations (PDEs) with multiscale coefficients. Inspired by the deep mixed residual method, we rewrite the second-order problem into a first-order system and employ multiple radial basis function neural networks (RBFNNs) to approximate unknown functions in the system. To aviod the overfitting due to the simplicity of RBFNN, an additional regularization is introduced in the loss function. Thus the loss function contains two parts: the $L_2$ loss for the residual of the first-order system and boundary conditions, and the $\ell_1$ regularization term for the weights of radial basis functions (RBFs). An algorithm for optimizing the specific loss function is introduced to accelerate the training process. The accuracy and effectiveness of the proposed method are demonstrated through a collection of multiscale problems with scale separation, discontinuity and multiple scales from one to three dimensions. Notably, the $\ell_1$ regularization can achieve the goal of representing the solution by fewer RBFs. As a consequence, the total number of RBFs scales like $\mathcal{O}(\varepsilon^{-n\tau})$, where $\varepsilon$ is the smallest scale, $n$ is the dimensionality, and $\tau$ is typically smaller than $1$. It is worth mentioning that the proposed method not only has the numerical convergence and thus provides a reliable numerical solution in three dimensions when a classical method is typically not affordable, but also outperforms most other available machine learning methods in terms of accuracy and robustness. 
%%%%
\end{abstract}

\date{September 11, 2020}
\maketitle
%\begin{keyword}
%\KWD Multiscale elliptic problem\sep Radial basis function neural network\sep $\ell_1$ regularization
%\end{keyword}

%\linenumbers

%% main text
 
\section{Introduction}
Multiscale phenomena are common in nature and mathematical modeling typically involves a small parameter $\varepsilon$ which characterizes the ratio between the smallest characteristic length over the size of interest \cite{weinan2011principles}. In some scenarios, the mechanical response of composite materials and flow property of porous media can be modeled by the elliptic problems with multiscale coefficients
\begin{equation}
	\begin{aligned}
		\left\{\begin{aligned}
			-\mbox{div} (a^{\varepsilon}(\textbf{x})\nabla u^{\varepsilon}(\textbf{x}))&=f(\textbf{x})\quad \textbf{x}\in \Omega,  \\
			u^{\varepsilon}(\textbf{x})&=g(\textbf{x})  \quad \textbf{x}\in \partial\Omega,
		\end{aligned}
		\right.
	\end{aligned}
	\label{eq1}
\end{equation}
where $\Omega\subset\mathbb{R}^n, n=1, 2, 3$ is a bounded domain, $\varepsilon$ $\ll 1$, and $f:\Omega\rightarrow\mathbb{R}$ is the source term.  

For small $\varepsilon$, classical numerical methods cannot solve \eqref{eq1} in an effective manner, where the number of unknowns scales like $\mathcal{O}(\varepsilon^{-n\tau})$ and $\tau>1$. Meanwhile, by the homogenization theory \cite{homogenization,mchomogenization}, as $\varepsilon\rightarrow 0$, an effective model can be obtained from \eqref{eq1} based on the periodicity assumption of $a^{\varepsilon}(\textbf{x})$.
Over the past a couple of decades, a large number of multiscale methods have been developed for multiscale elliptic problems in the literature, such as the multiscale finite element method (MsFEM) \cite{MsFEM0,MsFEM1}, wavelet homogenization techniques \cite{wavelet}, the heterogeneous multiscale method (HMM) \cite{HMM2,HMM3}, the metric-based upscaling \cite{owhadi2007metric}, and the local orthogonal decomposition (LOD) \cite{localization}. Typically, these methods solve multiple local problems in advance to extract the macroscopic information and then solve one global problem. The total computational cost is at least $\mathcal{O}(\varepsilon^{-n})$.

In recent years, a series of deep neural network methods for solving PDEs are derived. Deep Ritz method (DRM) \cite{DRM} employs the variational formulation as the loss function. Deep Galerkin method (DGM) \cite{DGM} and Physics-informed neural networks (PINNs) \cite{PINN} utilize the residual of the PDE as the loss function. Weak Adversarial Network (WAN) \cite{zang2020weak} employs the weak formulation to construct the loss function. Multilayer perception (MLP) and residual neural network (ResNet) are commonly used to approximate the solution. For multiscale problems, DNN-based methods are also developed. Multiscale deep neural network (MscaleDNN) \cite{MDNN1} approximates the solution by converting the original data to a low frequency space and \cite{MDNN2} improves the MscaleDNN algorithm by a smooth and localized activation function.

\begin{figure}[H]
	%\vspace{-1.0cm} 
	\centering
	\includegraphics[width=0.5\textwidth]{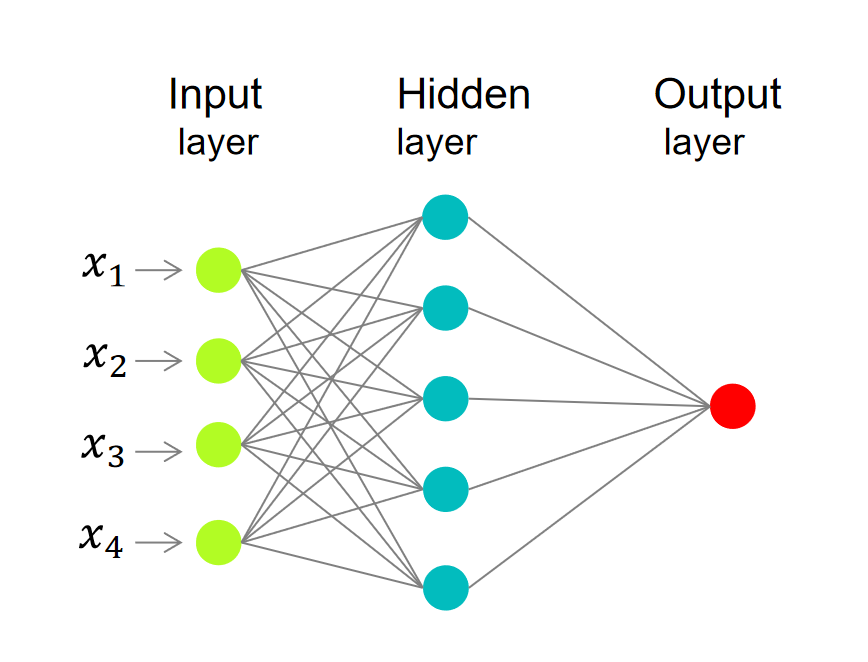}
	\caption{The network architecture of radial basis function neural network.}
	\label{RBFNN}
	%\vspace{-0.8cm}
\end{figure}

Radial basis function neural network (RBFNN) represents an attractive alternative to other neural network models \cite{RBFNNs}. One reason is that it forms a unified connection between function approximation, regularization, noisy interpolation, classification, and density estimation. It is also the case that training a RBFNN is faster than training MLP networks \cite{RBFNNs}. RBFNN is a simple three-layer feedforward neural network. The first layer is the input layer, representing the input features. The second layer is the hidden layer, composed of radial basis activation functions. The third layer is the output layer. Fig. \ref{RBFNN} shows the network architecture of RBFNN in four dimensions with five neurons.     Like most neural networks, RBFNN also has the universal approximation property \cite{park1991universal}. Meanwhile,
compared with deep neural networks, RBFNN avoids the tedious and lengthy calculation of back propagation between the input layer and the output layer, significantly accelerating the training process of the network. 

Often, optimizing the empirical risk leads to the overfitting problem owing to the simplicity of RBFNN, and causes a poor generalization capability. To tackle this issue, $\ell_1$ regularization techniques \cite{regular1,regular2}, become common components in many active fields, such as sparse representations of images \cite{chenzhaodi}, molecular surface \cite{lubenzhuo} and point cloud surface reconstruction \cite{wangdandan}. The basic idea is to add a regularization term to the
loss function to penalize over-complicated solutions. For $\ell_1$ regularization, a weighted $\ell_1$ norm of the parameter vector is added to the loss function, which penalizes the sum of the absolute values of the parameters.

In this work, we propose a sparse radial basis function neural network (SRBFNN) method for solving multiscale elliptic equations, where RBFNN is used to approximate the multiscale solution and the $\ell_1$ regularization term is added to the loss function to guarantee the sparsity of RBF representation. In the modeling stage, by rewriting the second-order equation into a first-order system as shown in the deep mixed residual method \cite{MIM}, we obtain an augmented problem in the sense that both the solution and the auxiliary variables related to first-order derivatives of the solution are unknown functions to be approximated.  An algorithm for optimizing the specific loss function is introduced to accelerate the training process. The accuracy and effectiveness of the proposed method are demonstrated through a collection of multiscale problems with scale separation, discontinuity and multiple scales from one to three dimensions. Notably, the $\ell_1$ regularization can achieve the goal of representing the solution by fewer RBFs. As a consequence, the total number of RBFs scales like $\mathcal{O}(\varepsilon^{-n\tau})$, where $\varepsilon$ is the smallest scale, $n$ is the dimensionality, and $\tau$ is typically smaller than $1$.

The paper is organized as follows. Section \ref{sec2} introduces the SRBFNN. An algorithm is proposed to optimize the loss function with the $\ell_1$ regularization during the training process in Section \ref{sec3}. Section \ref{sec4} provides numerical experiments from one to three dimensions to demonstrate the accuracy and effectiveness of the proposed method. Conclusions and discussions are drawn in Section \ref{sec5}.
 
\section{Sparse radial basis function neural network }
\label{sec2}
\subsection{Radial basis function }
The function $\Phi : \mathbb{R}^n \rightarrow \mathbb{R} $ is called a radial basis function, which can be passed through a one-dimensional (1D) function $\phi:[0, \infty) \rightarrow \mathbb{R}$, this means
\begin{eqnarray}
	\Phi(\textbf{x})=\phi(\Vert\textbf{x}\Vert)=\phi(r),\quad r:=\Vert\textbf{x}\Vert, \nonumber
\end{eqnarray}
where $\Vert\cdot\Vert$ represents the $L_2$ norm in $\mathbb{R}^n$.\\
Some common types of RBFs:
\begin{itemize}
	\item Gaussian function
	\begin{eqnarray}
		\Phi(r)=e^{-d^2r^2}; \nonumber
	\end{eqnarray}
	\item Multiquadric function (MQ)
	\begin{eqnarray}
		\Phi(r)=\sqrt{1+d^2r^2}; \nonumber
	\end{eqnarray}
	\item Inverse multiquadric function (IMQ)
	\begin{eqnarray}
		\Phi(r)=1/\sqrt{1+d^2r^2}, \nonumber
	\end{eqnarray}
\end{itemize}
where $d$ is the shape parameter that plays an important role in the approximation accuracy. As shown in \cite{davydov2011optimal}, a small shape parameter results in an ill-conditioned system, while the corresponding RBF method exhibits the high approximation accuracy. In contrast, a large shape parameter results in an well-conditioned system, but the approximation accuracy is low.

In this work, we use the ellipse Gaussian RBFs
\begin{eqnarray}	\Phi(\textbf{x};\textbf{c},\textbf{D})=\mbox{e}^{-\Vert\textbf{D}(\textbf{x}-\textbf{c})\Vert},
	\label{ellipserbf}
\end{eqnarray}
where $\textbf{D}=\mathrm{diag}(d_1,...,d_n)$ is the shape parameter and $\textbf{c}=(c_1,...,c_n)^T$ is the center. The ellipse Gaussian RBF has the good approximation property in a specific neighborhood of the center with different radius.

\subsection{Sparse radial basis function neural network}
The RBFNN can be expressed as
\begin{eqnarray}
	\textbf{NN}(\textbf{x};\textbf{w},\textbf{c},\textbf{D})=\sum\limits_{i=1}^{N}w_i\Phi_i(\textbf{x};\textbf{c}_i,\textbf{D}_i),
	\label{ellgaussian}
\end{eqnarray}
where $\Phi_i(\textbf{x};\textbf{c}_i,\textbf{D}_i), i=1,...,N$ are the ellipse Gaussian RBFs and $\textbf{w}=(w_1,..,w_N)^T$ is the corresponding weight.
For the ease of exposition, we denote the trainable parameters $\textbf{w},\textbf{c},\textbf{D}$ by $\boldsymbol{\theta}$, i.e.,
\begin{eqnarray}
	\boldsymbol{\theta} = \{\textbf{w}\in\mathbb{R}^N,\textbf{c}\in\mathbb{R}^{N\times n},\textbf{D}\in\mathbb{R}^{N\times n}|N \mbox{ is the number of basis, }n\mbox{ is the dimension of the input } \textbf{x}\}.\nonumber
\end{eqnarray}
The RBFNN whose loss function contains the $\ell_1$ regularization term for the weight vector \textbf{w} is called the sparse RBFNN (SRBFNN). 

\subsection{Loss function }
We consider multiscale elliptic equations of the form \eqref{eq1}. The total loss function for SRBFNN is defined as
\begin{eqnarray}	Loss(\textbf{x};\boldsymbol{\theta})=L_s(\textbf{x};\boldsymbol{\theta})+\lambda_{3}\Vert \textbf{w}^u\Vert_1,
\end{eqnarray}
where $L_s$ is the $L_2$ loss function for the first-order system, $\lambda_{3}$ is the penalty parameter, and $\textbf{w}^u$ is the weight of the RBFNN that approximates the solution. $\Vert\cdot\Vert_1$ denotes the $\ell_1$ norm.
Explicit forms of $L_s$ from one to three dimensions will be specified in \eqref{loss1d1}, \eqref{loss2d1}, and \eqref{loss3d1}.

\subsubsection{$L_2$ loss function in 1D}
Inspired by the MIM method \cite{MIM}, we rewrite \eqref{eq1} in 1D, using an auxiliary variable $p$, into a first-order system:
\begin{equation}
	\begin{aligned}
		\left\{\begin{aligned}
			p-a^\varepsilon u_x =0 \quad x\in\Omega,\\
			p_x+f =0\quad x\in\Omega,\\
			u=g\quad x\in\partial\Omega.
		\end{aligned}
		\right.
	\end{aligned}
	\label{system1d}
\end{equation}
Then we use two RBFNNs ($\textbf{NN}_1$, $\textbf{NN}_2$) to approximate $u$ and $p$ as
\begin{equation}
	\widetilde{u} = \textbf{NN}_1(x;\boldsymbol{\theta}^u), \quad\widetilde{p} =  \textbf{NN}_2(x;\boldsymbol{\theta}^p),
	\label{eq2.4.1}
\end{equation}
where $\boldsymbol{\theta}^u$ and $\boldsymbol{\theta}^p$ are the neural network parameters used to approximate $u$ and $p$, respectively. Then the $L_2$ loss function for the system \eqref{system1d} is defined by
\begin{equation}
	L_s(\textbf{x};\boldsymbol{\theta}^u,\boldsymbol{\theta}^p) =\Vert a^{\varepsilon}\widetilde{u}_x -\widetilde{p}\Vert^2_{2,\Omega}+\lambda_{1}\Vert \widetilde{p}_x +f\Vert^2_{2,\Omega} +\lambda_{2}\Vert\widetilde{u}-g\Vert^2_{2,\partial\Omega},
	\label{loss1d1}
\end{equation}
where $\lambda_{1}$ and $\lambda_{2}$ are penalty parameters given in advance. The three terms in \eqref{loss1d1} measure how well the approximate solution and the auxiliary variable satisfy the first-order system \eqref{system1d} and the boundary condition, respectively. Due to the simple form of the RBFNN defined in \eqref{ellgaussian}, the corresponding derivatives $\widetilde{u}_x$, $\widetilde{p}_x$ can be easily calculated.

\subsubsection{$L_2$ loss function in 2D}
Analogously, we rewrite \eqref{eq1} in 2D, using auxiliary variables $p,q$, as
\begin{equation}
	\begin{aligned}
		\left\{\begin{aligned}
			p-a^{\varepsilon}\cdot \frac{\partial u}{\partial x}=0\quad \textbf{x}\in\Omega,\\
			q-a^{\varepsilon}\cdot \frac{\partial u}{\partial y}=0\quad \textbf{x}\in\Omega,\\
			\frac{\partial p}{\partial x} + \frac{\partial q}{\partial y} +f=0 \quad \textbf{x}\in\Omega,\\
			u=g \quad \textbf{x}\in\partial\Omega.
		\end{aligned}
		\right.
	\end{aligned}
	\label{eq6}
\end{equation}
Then we use three RBFNNs ($\textbf{NN}_1$, $\textbf{NN}_2$, $\textbf{NN}_3$) to approximate $u, p, q$ as
\begin{equation}
	\widetilde{u} = \textbf{NN}_1(\textbf{x};\boldsymbol{\theta}^u),\quad\widetilde{p} =  \textbf{NN}_2(\textbf{x};\boldsymbol{\theta}^p),\quad\widetilde{q} =  \textbf{NN}_3(\textbf{x};\boldsymbol{\theta}^q).
	\label{eq2.4.2}
\end{equation}
The $L_2$ loss function for the system \eqref{eq6} is defined as
\begin{eqnarray}
	L_s(\textbf{x};\boldsymbol{\theta}^u,\boldsymbol{\theta}^p,\boldsymbol{\theta}^q) =\Vert a^{\varepsilon}\cdot\frac{\partial\widetilde{u}}{\partial x} -\widetilde{p}\Vert^2_{2,\Omega}+\Vert a^{\varepsilon}\cdot\frac{\partial \widetilde{u}}{\partial y}-\widetilde{q}\Vert^2_{2,\Omega}+\lambda_{1}\Vert \frac{\partial \widetilde{p}}{\partial x}+\frac{\partial \widetilde{q}}{\partial y}+f\Vert^2_{2,\Omega} +\lambda_{2}\Vert\widetilde{u}-g\Vert^2_{2,\partial\Omega}.
	\label{loss2d1}
\end{eqnarray}

\subsubsection{$L_2$ loss function in 3D}
We rewrite \eqref{eq1} in 3D, using auxiliary variables $p,q,r$, as
\begin{equation}
	\begin{aligned}
		\left\{\begin{aligned}
			p-a^{\varepsilon}\cdot \frac{\partial u}{\partial x}=0 \quad \textbf{x}\in\Omega,\\
			q-a^{\varepsilon}\cdot \frac{\partial u}{\partial y}=0\quad \textbf{x}\in\Omega,\\
			r-a^{\varepsilon}\cdot \frac{\partial u}{\partial z}=0\quad \textbf{x}\in\Omega,\\
			\frac{\partial p}{\partial x} + \frac{\partial q}{\partial y}+ \frac{\partial r}{\partial z} +f=0\quad \textbf{x}\in\Omega,\\
			u=g \quad \textbf{x}\in\partial\Omega.
		\end{aligned}
		\right.
	\end{aligned}
	\label{eq3d}
\end{equation}
Then we use four RBFNNs ($\textbf{NN}_1$, $\textbf{NN}_2$, $\textbf{NN}_3$, $\textbf{NN}_4$) to approximate $u, p, q,r$ as
\begin{equation}
	\widetilde{u} = \textbf{NN}_1(\textbf{x};\boldsymbol{\theta}^u),\quad\widetilde{p} =  \textbf{NN}_2(\textbf{x};\boldsymbol{\theta}^p),\quad\widetilde{q} =  \textbf{NN}_3(\textbf{x};\boldsymbol{\theta}^q),\quad\widetilde{r} =  \textbf{NN}_4(\textbf{x};\boldsymbol{\theta}^r).
	\label{eq2.4.3}
\end{equation}
The $L_2$ loss function for the system \eqref{eq3d} is 
\begin{equation}
	\begin{split}
		L_s(\textbf{x};\boldsymbol{\theta}^u,\boldsymbol{\theta}^p,\boldsymbol{\theta}^q,\boldsymbol{\theta}^r) = &\Vert a^{\varepsilon}\cdot\frac{\partial\widetilde{u}}{\partial x} -\widetilde{p}\Vert^2_{2,\Omega}+\Vert a^{\varepsilon}\cdot\frac{\partial \widetilde{u}}{\partial y}-\widetilde{q}\Vert^2_{2,\Omega}+\Vert a^{\varepsilon}\cdot \frac{\partial\widetilde{u}}{\partial z} -\widetilde{r}\Vert^2_{2,\Omega}+\\
		&\lambda_{1}\Vert \frac{\partial \widetilde{p}}{\partial x}+\frac{\partial \widetilde{q}}{\partial y}+\frac{\partial \widetilde{r}}{\partial z}+f\Vert^2_{2,\Omega} +\lambda_{2}\Vert\widetilde{u}-g\Vert^2_{2,\partial\Omega}.
		\label{loss3d1}
	\end{split}
\end{equation}

\subsection{Network initialization}
All RBFNNs use the same initialization strategy for the parameters $\textbf{w}$, $\textbf{c}$, $\textbf{D}$. The center $\textbf{c}$ is initialized with a uniform distribution over $\Omega$. For the weight $\textbf{w}$, we initialize them with a uniform distribution U(0,1). The initialization of the shape parameter $\textbf{D}$ is extremely important for the network convergence in multiscale problems. According to the homogenization theory \cite{homogenization}, the solution $u^{\varepsilon}$ to \eqref{eq1} admits the following asymptotic expansion:
\begin{eqnarray}
	u^\varepsilon(\textbf{x}) = u_0(\textbf{x})+\varepsilon u_1(\textbf{x},\frac{\textbf{x}}{\varepsilon})+\varepsilon^2u_2(\textbf{x},\frac{\textbf{x}}{\varepsilon})+\cdots,
\end{eqnarray}
where $u_0$ is the homogenized solution, and $u_1$, $u_2$,... are high-order terms and depend on the derivatives of $u_0$. Therefore, to resolve the smallest scale of the solution, we initialize $\textbf{D}$ with uniform distribution U(0, $\frac{1}{\varepsilon}$).

\section{Training process of SRBFNN}
\label{sec3}
In the training process of SRBFNN, we aim to find a sparse weight vector, $\textbf{w}^u$, while $\textbf{NN}_i, i=1, ..., n+1$ are also good fits of the solution and the auxiliary variables. Therefore, an algorithm to optimize the loss function with $\ell_1$ regularization during the training phase is proposed in Algorithm \ref{algorithm1}. 

\begin{algorithm}[H]
	\SetKwInOut{Input}{Input}
	\SetKwInOut{Output}{Output}
	\caption{Training process of the SRBFNN}\label{algorithm1}
	\Input{The coefficient $a^\varepsilon$, the source term $f$, and the boundary condition $g$ in \eqref{eq1}.}
	\Output{The list of parameters for $\textbf{NN}_i, i=1, ..., n+1$.}
	\begin{enumerate}[$\textbf{Step}$ 1.] 
		\item initialize $\textbf{NN}_i$, $i=1,...,n+1$, $n$ is the dimension\;
		\item generate two sets of points uniformly distributed over $\Omega$ and $\partial\Omega$: $\{\textbf{x}_m\}_{m=1}^{M_r}$ in $\Omega$ and  $\{\widehat{\textbf{x}}_m\}_{m=1}^{M_b}$ on $\partial\Omega$, $M_r$ and $M_b$ are the number of sample points\;
		\item initialize some hyperparameters\;
		\begin{enumerate}[$\textbf{Step}$ 3.1.]
			\item set the penalty coefficients of the loss functions, i.e., $\lambda_{1}$, $\lambda_{2}$ and $\lambda_{3}$\;
			\item set the maximum iteration $MaxNiter$ and the number of sparse optimization iterations $SparseNiter$. Set $Check_{iter}$ for deleting the neglectable weight every $Check_{iter}$ iterations\;
			\item initialize the threshold, $thres$=0.0, for adding $\ell_1$ regularization and the recorded loss, $L^{rec}$=0.0, for recording latest $L_2$ loss per $Check_{iter}$ iterations. Initialize the learning rate\;  
			\item set the tolerance, $tol_1$, to determine whether $L_2$ loss is convergent and the tolerance, $tol_2$, for neglectable weight in $\textbf{NN}_1$. Set $batchsize$, $Niter=0$ and $j=0$\;
		\end{enumerate}
		\item	start following training process\;
	\end{enumerate}
	\While{$Niter\leq MaxNiter$}{
		$Niter+=1$\; 
		\For{$i=0,i\leq M_r, i+=batchsize$}{
			$\textbf{Step}$ 4.1.
			calculate the $L_2$ loss function, $L_s$ (in \ref{loss1d1}, \ref{loss2d1}, \ref{loss3d1}), with $batchsize$ sample points\;
			$\textbf{Step}$ 4.2. \uIf{$L_s>thres$ {\bf or} $Niter>SparseNiter$}{$\lambda_{3}=0$\;}\Else{
			$\lambda_{3}>0$\;}
			$\textbf{Step}$ 4.3. update the network parameters by ADAM optimizer\;
		}
		\If{$Niter\%Check_{iter}==0$}{
			$\textbf{Step}$ 4.4. \If{$|L^{rec}-L_s|<tol_1$ {\bf and} $j==0$}{
				set $thres=L_s+tol_1$\;
				$j+=1$\;}
			$\textbf{Step}$ 4.5. \If{$j>0$ {\bf and} $Niter<SparseNiter$}{ delete neglectable basis function, whose coefficient is close to zero, i.e., $|w_i| < tol_2$\;}
			$\textbf{Step}$ 4.6. record the $L_2$ loss as the $L^{rec}$ every $Check_{iter}$ iterations, i.e., $L^{rec}=L_s$ \;
		}
	}
	
\end{algorithm}
Step 1 initializes the parameters of SRBFNN. Step 2 generates two sets of training points $\{\textbf{x}_m\}_{m=1}^{M_r}$ and $\{\widehat{\textbf{x}}_m\}_{m=1}^{M_b}$ uniformly distributed over $\Omega$ and $\partial\Omega$, respectively.
Step 3.1 sets the penalty coefficients $\lambda_{1}$, $\lambda_{2}$, $\lambda_{3}$. Step 3.2 sets the number of
maximum iterations and the number of sparse optimization
iterations, i.e., the max number of iterations with the loss function having the $\ell_1$ regularization term. Step 3.3 sets a threshold to the $L_2$ loss, which denotes the sign for adding the $\ell_1$ regularization, i.e., setting $\lambda_{3}>0$. Step 3.4 sets a tolerance for the stability condition of $L_2$ loss.
Step 4.1 calculates the $L_2$ loss with $batchsize$ sample points. When the $L_s$ is lower than $thres$ and $Niter <SparseNiter$, we add the $\ell_1$ regularization in step 4.2. Step 4.3 selects the ADAM optimizer \cite{Adam} to update the parameters and its pipeline is as follows
\begin{equation}
	\begin{aligned}
		m_k &= \beta_1\cdot m_{k-1}+(1-\beta_1)\cdot\nabla Loss_{k},\nonumber\\
		v_k &= \beta_2\cdot v_{k-1}+(1-\beta_2)\cdot\nabla (Loss_{k})^2,\nonumber\\
		\boldsymbol{\theta}_k &=\boldsymbol{\theta}_{k-1}-lr\times \frac{\frac{m_k}{1-\beta_1^k}}{\sqrt{\frac{v_k}{1-\beta_2^k}}+\epsilon},\nonumber
	\end{aligned}
\end{equation}
where $lr$ is the learning rate. $\beta_1$, $\beta_2$ and $\epsilon$ are set as default values ($\beta_1$ = 0.9, $\beta_2$ = 0.999 and $\epsilon$ = $10^{-8}$), $m_k$ is the $k$th biased
first moment estimate ($m_0$ = 0). $v_k$ is the $k$th biased second raw estimate ($v_0$ = 0).
Step 4.4 determines whether the difference between the $L_2$ loss at current iteration and the recorded loss is less than the tolerance, $tol_1$. If it does, $L_2$ loss is convergent and we let $thres$ equal to the $L_2$ loss of the current iteration plus $tol_1$. Step 4.5 determines whether $L_2$ loss is convergent and $Niter<SparseNiter$. If it does, we delete the neglectable basis function in $\boldsymbol{\theta}^u$, whose coefficient is close to 0. Step 4.6 records the value of $L_2$ loss per $Check_{iter}$ iterations. 

%----------Sparse radial basis function-------------------%

%----------Numerical examples for one dimension----------------%
\section{Numerical experiments }\label{sec4}
In this section, we present several numerical examples to illustrate the accuracy and effectiveness of our method. The coefficients are periodic functions in Example 1 and Example 5, two-scale problems with scale separation in Example 2, Example 3, and Example 6, a discontinuous periodic function in Example 4, a multiscale case in Example 7, respectively. Example 8 is a three dimensional example with the periodic coefficient.

We shall use the following three relative errors to measure the approximation accuracy.  
\begin{eqnarray} 
	err_2=\frac{\Vert u^S - u^F\Vert_{L_2(\Omega)}}{\Vert u^F\Vert_{L_2(\Omega)}}, \quad
	err_{\infty}=\frac{\Vert u^S - u^F\Vert_{L_\infty(\Omega)}}{\Vert u^F\Vert_{L_\infty(\Omega)}}, \quad 
	err_{H_1}=\frac{\Vert u^S - u^F\Vert_{H_1(\Omega)}}{\Vert u^F\Vert_{H_1(\Omega)}},
	\label{relative}
\end{eqnarray}
where $u^S$ is the approximate solution of SRBFNN and $u^F$ is the reference solution of finite difference method (FDM) with a very fine mesh size. The reference solution is not available in 3D since FDM is too expensive to solve a multiscale problem.
$\Vert\cdot\Vert_{L_2(\Omega)}$ is the $L_2$ norm, $\Vert\cdot\Vert_{L_\infty(\Omega)}$ is the $L_\infty$ norm and $\Vert\cdot\Vert_{H_1(\Omega)}$ is defined as
\begin{equation}
	\Vert u^S\Vert_{H_1(\Omega)} = \sqrt{\Vert u^S \Vert^2_{L_2(\Omega)}+\Vert  \nabla u^S \Vert^2_{L_2(\Omega)}}.
	\label{error}
\end{equation}
where $\nabla u^S$ is calculated by the auxiliary variables, i.e., $p/a^{\varepsilon}$ and $q/a^{\varepsilon}$. 

\subsection{Numerical experiments in 1D}
\textbf{Setting:}\label{1dset}
In all one-dimensional cases, unless otherwise specified, the initial learning rate is set to be 0.1 and reduced by 1/10 every 300 iterations until it is less than $10^{-5}$. We sample 10000 training points in the domain with a uniform distribution. The mesh size of FDM is set to 0.0001. $batchsize$ is 2048, $MaxNiter$ is 3000, $SparseNiter$ is 2000, $tol_1$ is $10^{-3}$, $tol_2$ is $10^{-5}$, $Check_{iter}$ is 100. The hyperparameters $\lambda_{1}$, $\lambda_{2}$ and $\lambda_{3}$ are initialized to 1.0, 100.0 and 0.001, respectively. 

We use six different values, $\varepsilon=0.5,0.1,0.05,0.01,0.005,0.002$, to check the accuracy of our method. The initial number of RBFs is set to 100, 200, 300, 500, 1000, 1500, respectively. The computational domain is the unit interval. $f(x)$ and $g(x)$ in \eqref{eq1} are set to 1.0 in all 1D cases. All experiments are conducted on the GPU device RTX 2080Ti.

\paragraph{Example 1}\label{example1}
Consider the case that $a^{\varepsilon}(x)$ is a periodic function
\begin{eqnarray}
	a^{\varepsilon}(x)=2+\sin(\frac{2\pi x}{\varepsilon}).
	\label{exam1}
\end{eqnarray}
Fig. \ref{figexam1_1} plots the solution and its derivative obtained by SRBFNN and FDM when $\varepsilon=0.01$, 0.002, respectively. One can see that SRBFNN has a good approximation for both the solution and the derivative.
\begin{small}
	\begin{figure}[H]
		\centering
		\subfigure[$\varepsilon=0.01$]{ 
			\includegraphics[align=c,width=0.4\columnwidth]{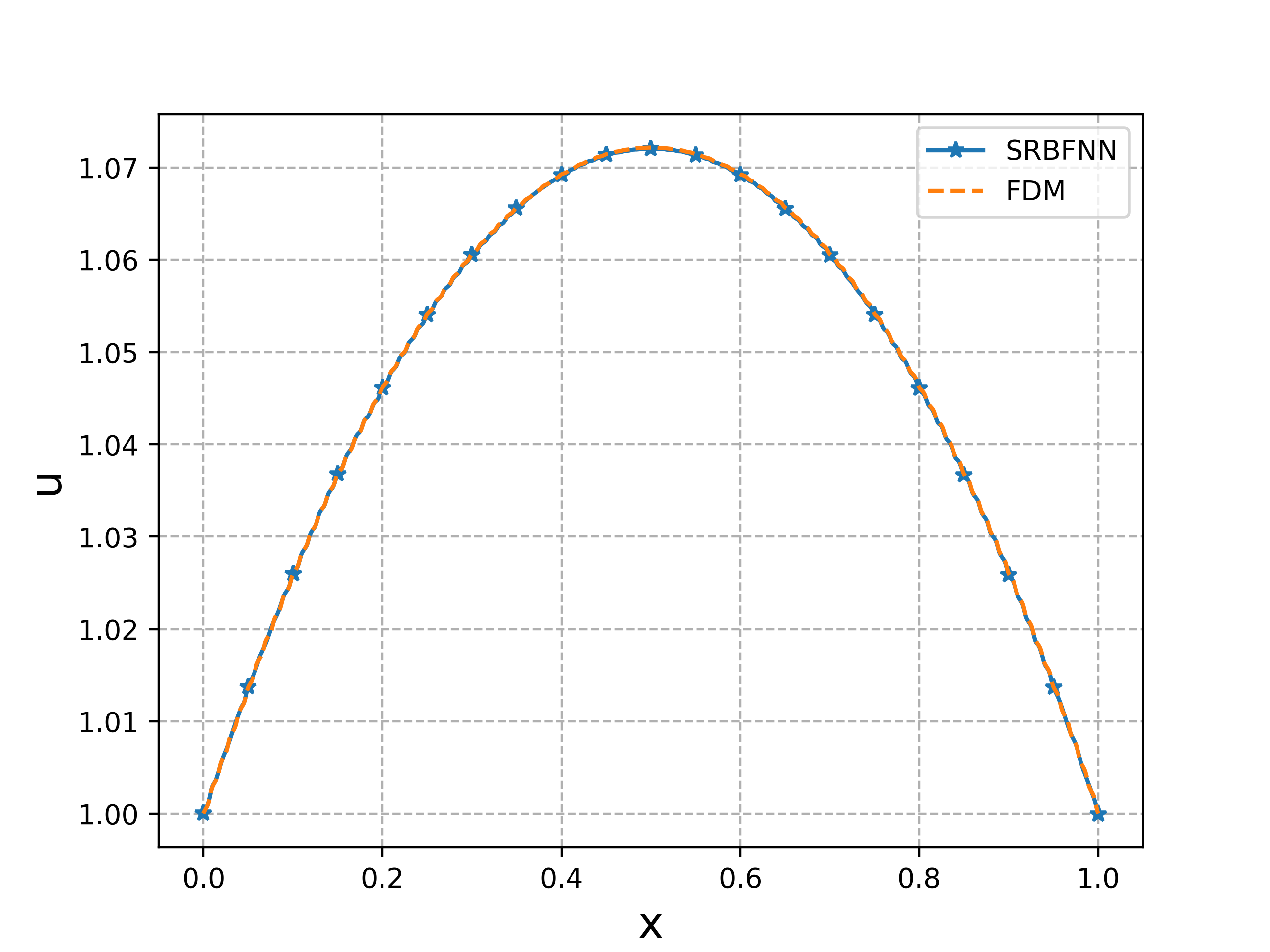}
		}
		\hfill 
		\subfigure[$\varepsilon=0.01$]{ 
			\includegraphics[align=c,width=0.4\columnwidth]{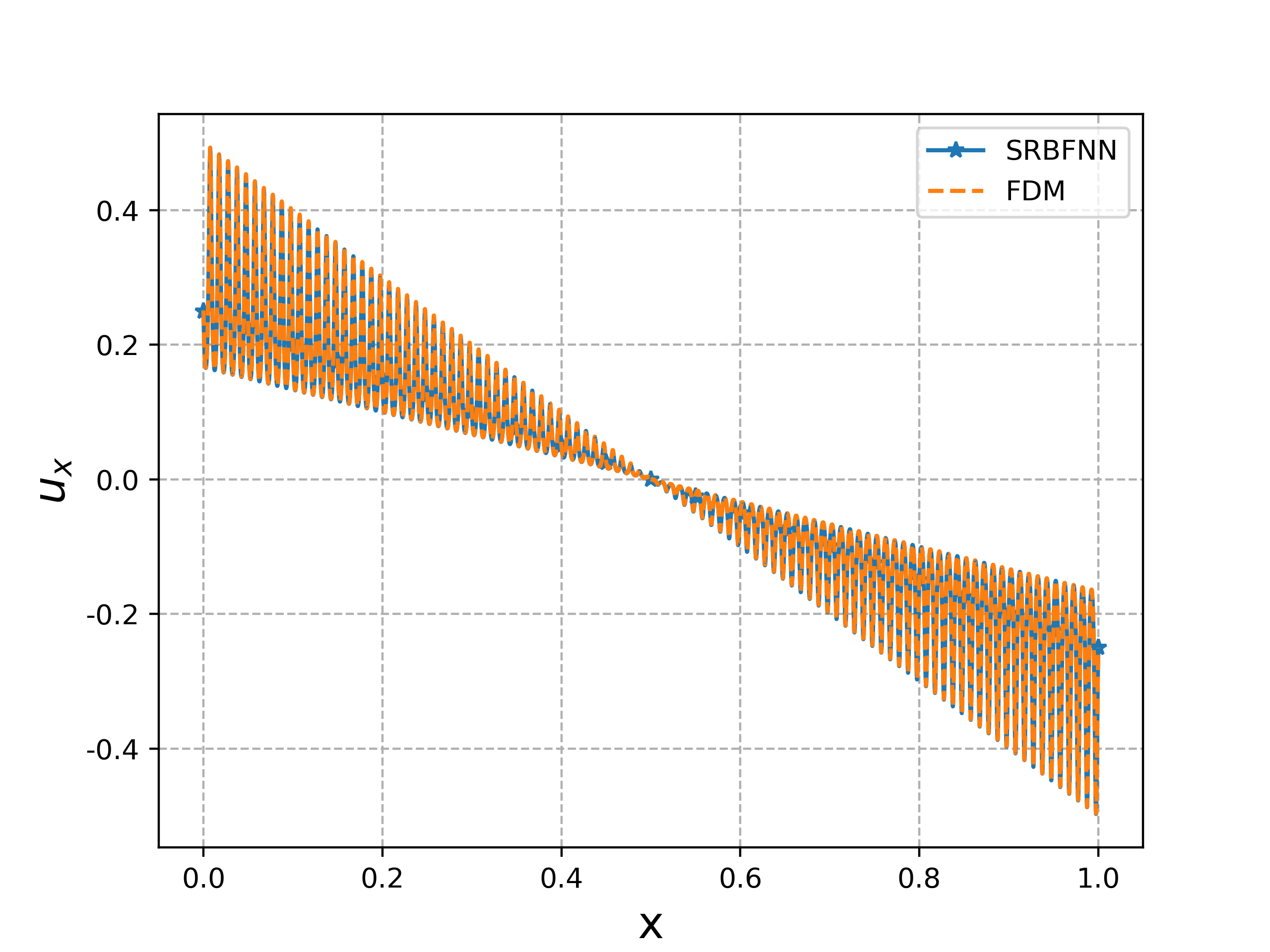}
		}
		\subfigure[$\varepsilon=0.002$]{ 
			\includegraphics[align=c,width=0.4\columnwidth]{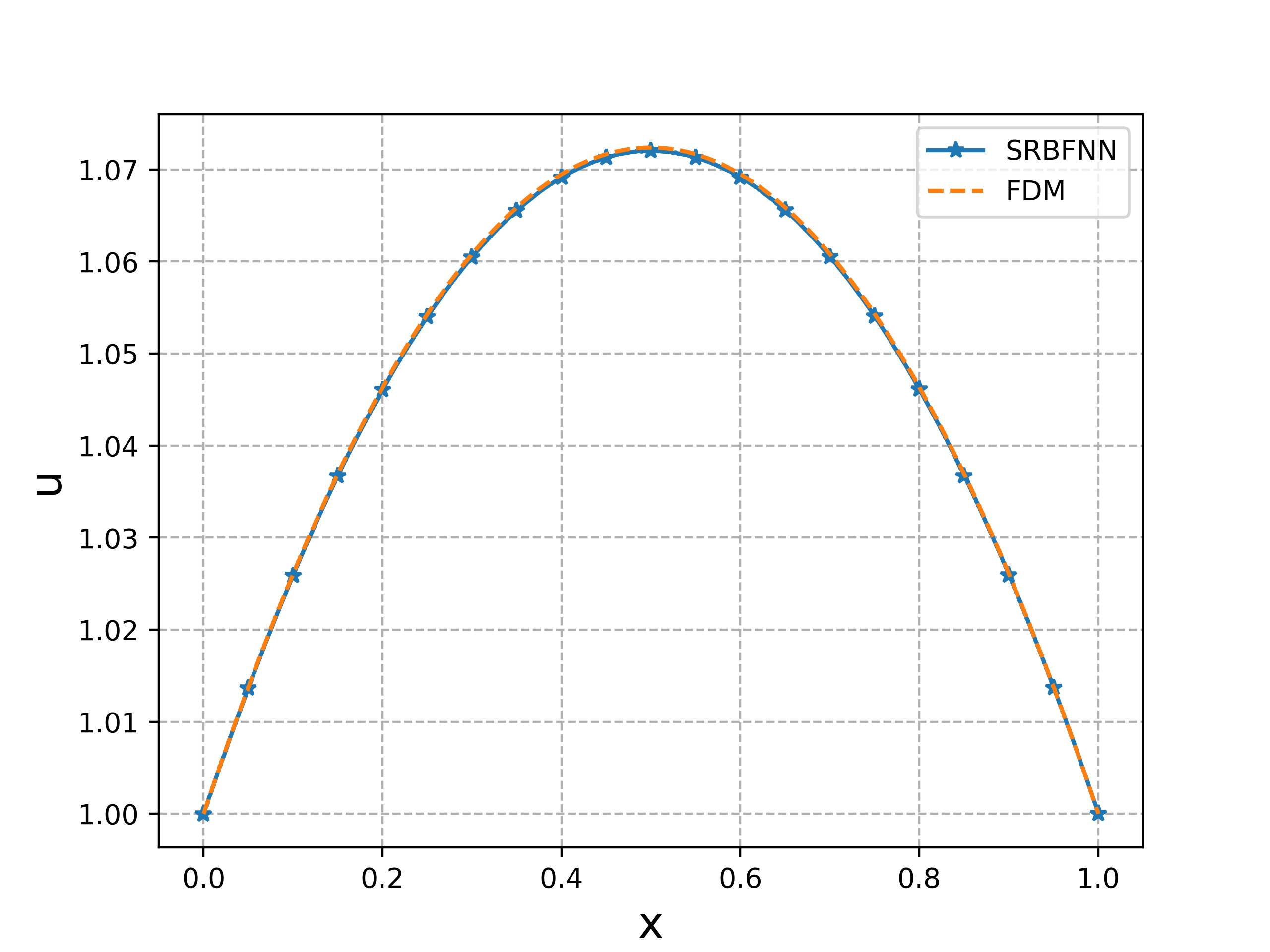}
		}
		\hfill 
		\subfigure[$\varepsilon=0.002$]{ 
			\includegraphics[align=c,width=0.4\columnwidth]{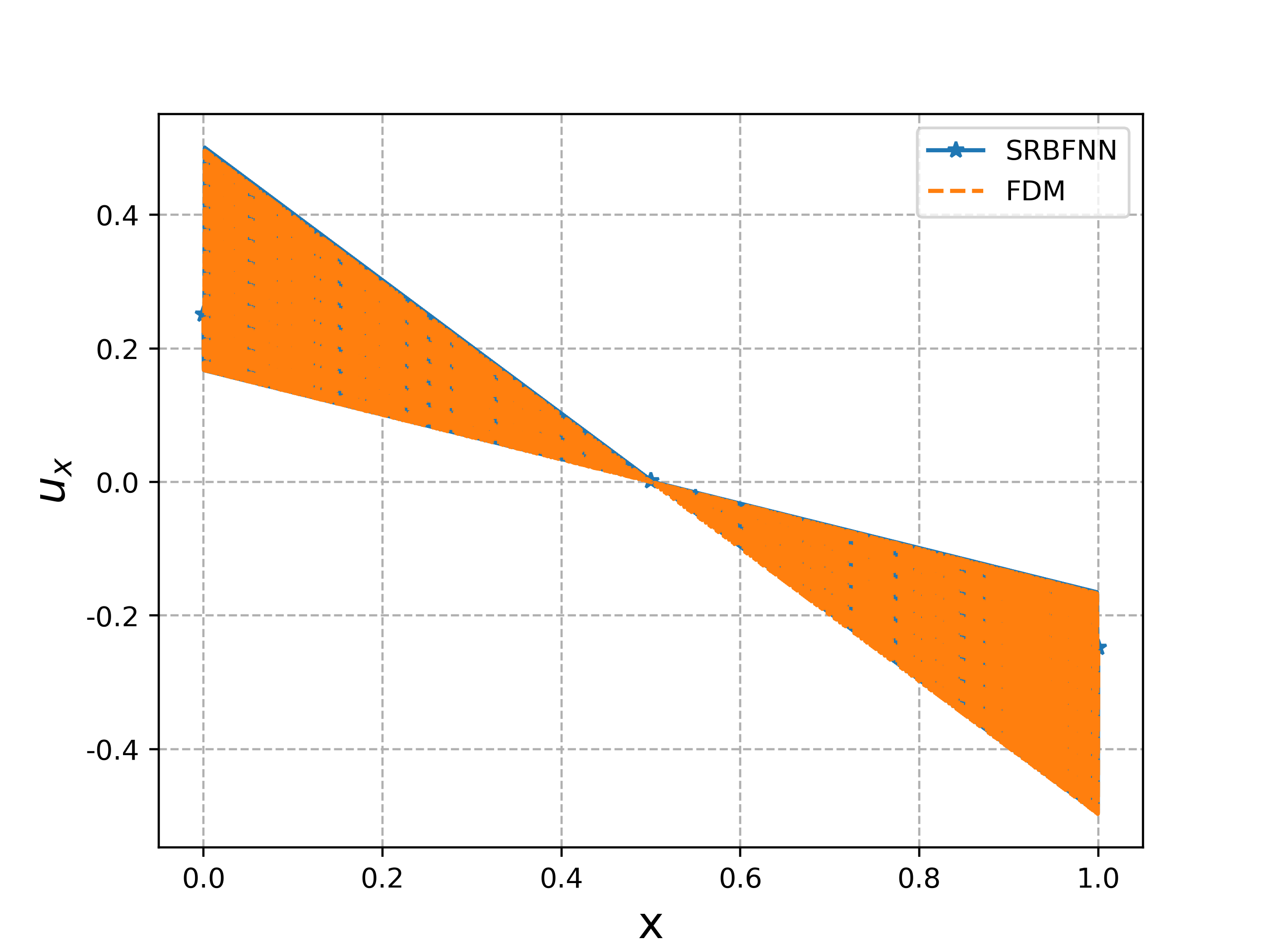}
		}
		\caption{Numerical solution and its derivative obtained by SRBFNN and FDM when $\varepsilon$=0.01, 0.002 for Example 1. }
		\label{figexam1_1}
	\end{figure}
\end{small}

Table \ref{tabexam1} shows three relative errors and the number of basis functions in the final solution from SRBFNN after the training process.
\begin{table}[H]
	\centering
	\begin{tabular}{ccccc}
		\hline
		$\varepsilon$ & $N$ & $err_2$ & $err_\infty$ & $err_{H_1}$  \\ \hline	
		\rowcolor{gray!40} 
		0.5   &  17 & 3.100e-5  & 8.552e-5 & 1.987e-4 \\
		0.1   &  38 & 6.718e-5  & 2.348e-4 & 8.560e-4 \\
		\rowcolor{gray!40}			
		0.05  &  66 & 8.311e-5  & 3.291e-4 & 7.847e-4 \\
		0.01  &  186& 1.086e-4  & 4.052e-4 & 8.433e-4 \\
		\rowcolor{gray!40}			
		0.005 &  367& 7.481e-5  & 4.107e-4 & 1.017e-3 \\
		0.002 &  750& 2.541e-4  & 3.940e-4 & 1.899e-3 \\\hline
	\end{tabular}
	\caption{Results of SRBFNN for Example 1. The second column records the number of basis functions in the final solution and the last three columns show the relative $L_2$,  $L_{\infty}$ and $H_1$ errors.}
	\label{tabexam1}
\end{table}

From Table \ref{tabexam1}, one general observation is that the number of RBFs needed in the final solution increases as the $\varepsilon$ decreases. In addition, the relative $L_2$ and $L_{\infty}$ errors from SRBFNN are of $10^{-4}$ for all scales and the relative $H_1$ error is about $10^{-3}$.

\paragraph{Example 2}
In this case, we consider the case that $a^\varepsilon(x)$ is a two-scale function with scale separation
\begin{eqnarray}
	a^{\varepsilon}(x)=2+\sin(\frac{2\pi x}{\varepsilon})\cos(2\pi x).
	\label{exam2}
\end{eqnarray}
We conduct a detailed comparison between SRBFNN and some deep neural network methods, including PINN \cite{PINN}, DGM \cite{DGM}, DRM \cite{DRM} and MscaleDNN \cite{MDNN2}. The training and test points are the same as those used in SRBFNN. PINN has five linear layers. DGM and DRM contain two residual blocks. The scale vector for MscaleDNN is [1, 2,..., $\frac{1}{\varepsilon}$]. The hidden layers are 16, 16, 32, 32, 64, 64 when $\varepsilon=$0.5, 0.1, 0.05, 0.01, 0.005, 0.002 for PINN, DGM and DRM. The hidden layers of MscaleDNN are (1000, 200, 150, 150, 100, 50, 50) for all $\varepsilon$. $tanh$ is used as the activation function for DRM, DGM and PINN. MscaleDNN employs a smooth and localized activation function $s2ReLU$, 
\begin{eqnarray}
	s2ReLU(x) = \sin(2\pi x)\ast\mbox{ReLU}(x)\ast\mbox{ReLU}(1-x). \nonumber
\end{eqnarray}  

Table \ref{tabexam2} shows the results of relative $L_2$ and $H_1$ errors for 5 methods at different scales.
\begin{table}[H]
	\centering
	\resizebox{.9\textwidth}{!}{
		\begin{tabular}{c|ccccc|ccccc}
			\toprule
			& &\multicolumn{3}{c}{$err_2$}& &&\multicolumn{3}{c}{$err_{H_1}$}  &\\ 
			$\varepsilon$ &SRBFNN&MscaleDNN&DRM&DGM&PINN&    SRBFNN&MscaleDNN&DRM&DGM&PINN\\ \hline	
			\rowcolor{gray!40}
			0.5  &6.602e-6&2.531e-4  &4.404e-2&1.697e-3&1.957e-3 & 1.977e-4& 3.983e-3 &1.437e-1&2.192e-2&2.122e-2\\
			0.1  &3.598e-5& 3.235e-4&5.237e-2&6.200e-3&6.168e-3 & 2.058e-4&1.111e-2&1.746e-1&4.876e-2&4.869e-2\\
			\rowcolor{gray!40}
			0.05 &4.774e-5&2.743e-3&2.874e-2&4.466e-3&6.262e-3 & 3.368e-4&4.177e-2&1.020e-1&4.895e-2&4.923e-2\\
			0.01 &2.016e-4&3.525e-3&2.874e-2&5.414e-3&5.549e-3 & 2.801e-3&4.694e-2&1.021e-1&4.924e-2&4.937e-2\\
			\rowcolor{gray!40}
			0.005&7.424e-5&3.695e-3&2.656e-2&7.424e-3&6.347e-3 & 1.282e-3&4.726e-2&9.597e-2&4.968e-2&4.935e-2\\
			0.002&1.654e-4&3.424e-3&2.661e-2&6.077e-3&5.511e-3 & 2.125e-3&4.711e-2&9.598e-2&4.914e-2&4.912e-2\\
			\bottomrule
		\end{tabular}
	}
	\caption{Comparison of SRBFNN, MscaleDNN, DRM, DGM, and PINN in terms of relative $L_2$ and $H_1$ errors for Example 2. The penalty coefficients of the boundary condition are 20.0, 20.0, 200.0 and 200.0 for PINN, DGM , DRM and MscaleDNN, respectively. The maximum number of iterations is 3000 for PINN, DGM, DRM and MscaleDNN.}
	\label{tabexam2}
\end{table}			
It is observed from Table \ref{tabexam2} that more accurate results can be obtained by SRBFNN and MscaleDNN. This may be explained by the fact that SRBFNN and MscaleDNN use the activation functions with good local approximation capability, which is an essential factor for multiscale problems. Next, we compare the number of parameters in final solution for all methods and the average time per iteration, i.e., $\frac{\mbox{Total running time}}{\mbox{Iterations}}$ in Table \ref{tabexam2_1}. It is found that SRBFNN contains fewer parameters and costs less time per iteration than other methods. This can be explained by the fact that SRBFNN has a simple network structure.

\begin{table}[H]
	\centering
	\resizebox{.9\textwidth}{!}{
		\begin{tabular}{c|ccccc|ccccc}
			\toprule
			&\multicolumn{5}{c}{Parameters }&\multicolumn{5}{c}{Average time per iteration (Seconds)} \\ 
			$\varepsilon$&SRBFNN&MscaleDNN&DRM&DGM&PINN & SRBFNN&MscaleDNN&DRM&DGM&PINN\\ \hline	
			\rowcolor{gray!40}		
			0.5  &351&277751 &1137 &1137 &865&0.122&0.714  &0.143  &0.877  &0.792  \\
			0.1  &699&277751 &1137 &1137 &865&0.140&0.746  &0.142  &0.872  &0.754  \\
			\rowcolor{gray!40}
			0.05 &1062&277751&4321&4321&3265 &0.147&0.764  &0.147  &0.843 &0.796  \\
			0.01 &1977&277751&4321&4321&3265 &0.144&0.793  &0.139  &0.861 &0.777  \\
			\rowcolor{gray!40}
			0.005&4107&277751&16833&16833&12673&0.149&0.825 &0.171  &0.799 &0.837  \\
			0.002&6474&277751&16833&16833&12673&0.150&0.734  &0.176  &0.841  &0.835  \\ \bottomrule
	\end{tabular}}
	\caption{Comparison of SRBFNN, MscaleDNN, DRM, DGM, and PINN in terms of the number of parameters and the averge time per iteration measured in seconds for Example 2.}
	\label{tabexam2_1}
\end{table}

\begin{figure}[H]
	\centering
	\subfigure[Relative error]{ 
		\includegraphics[align=c,width=0.45\columnwidth]{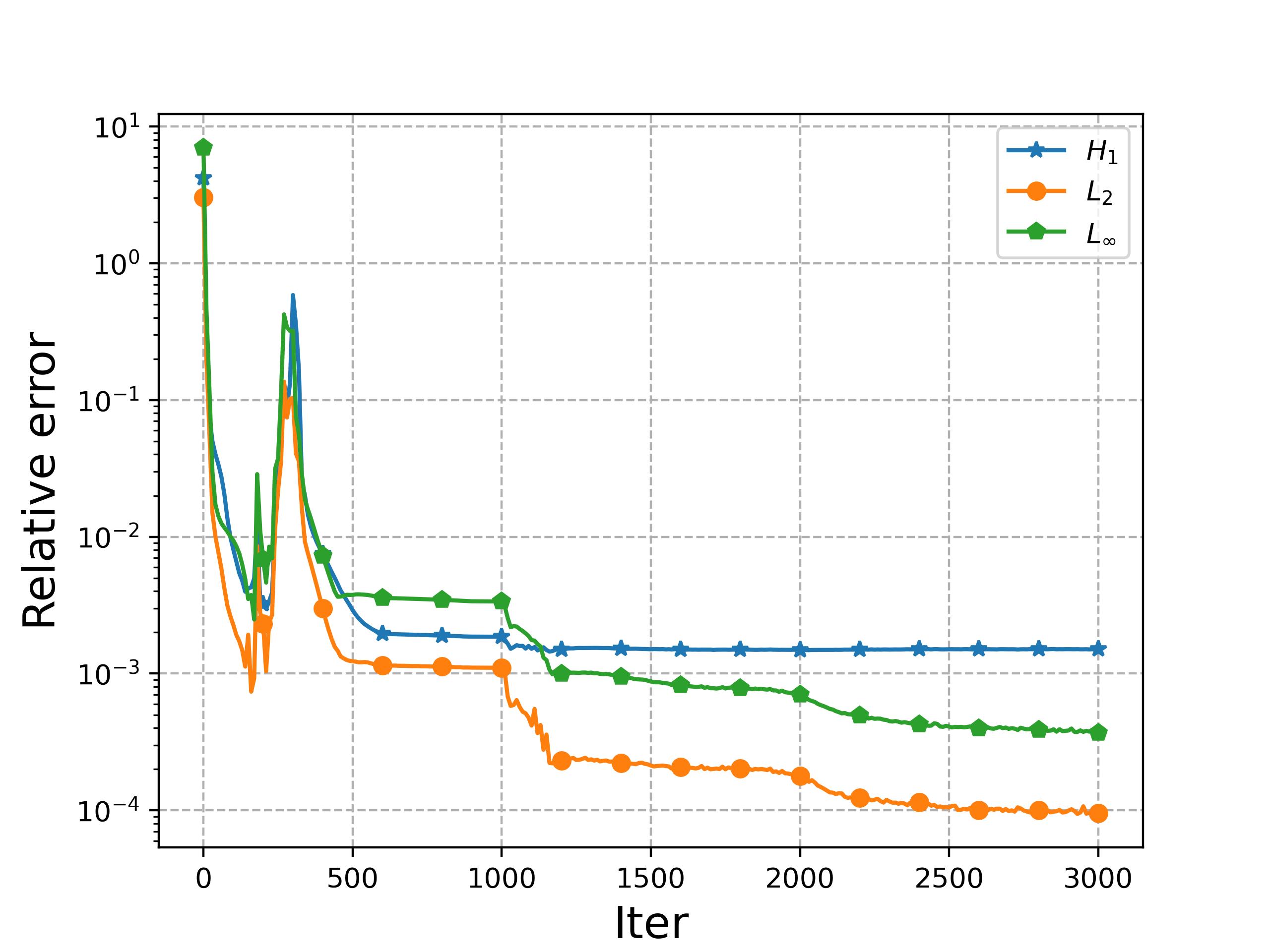}
	}
	\hfill %   ???   ?   
	\subfigure[The number of RBFs]{ 
		\includegraphics[align=c,width=0.45\columnwidth]{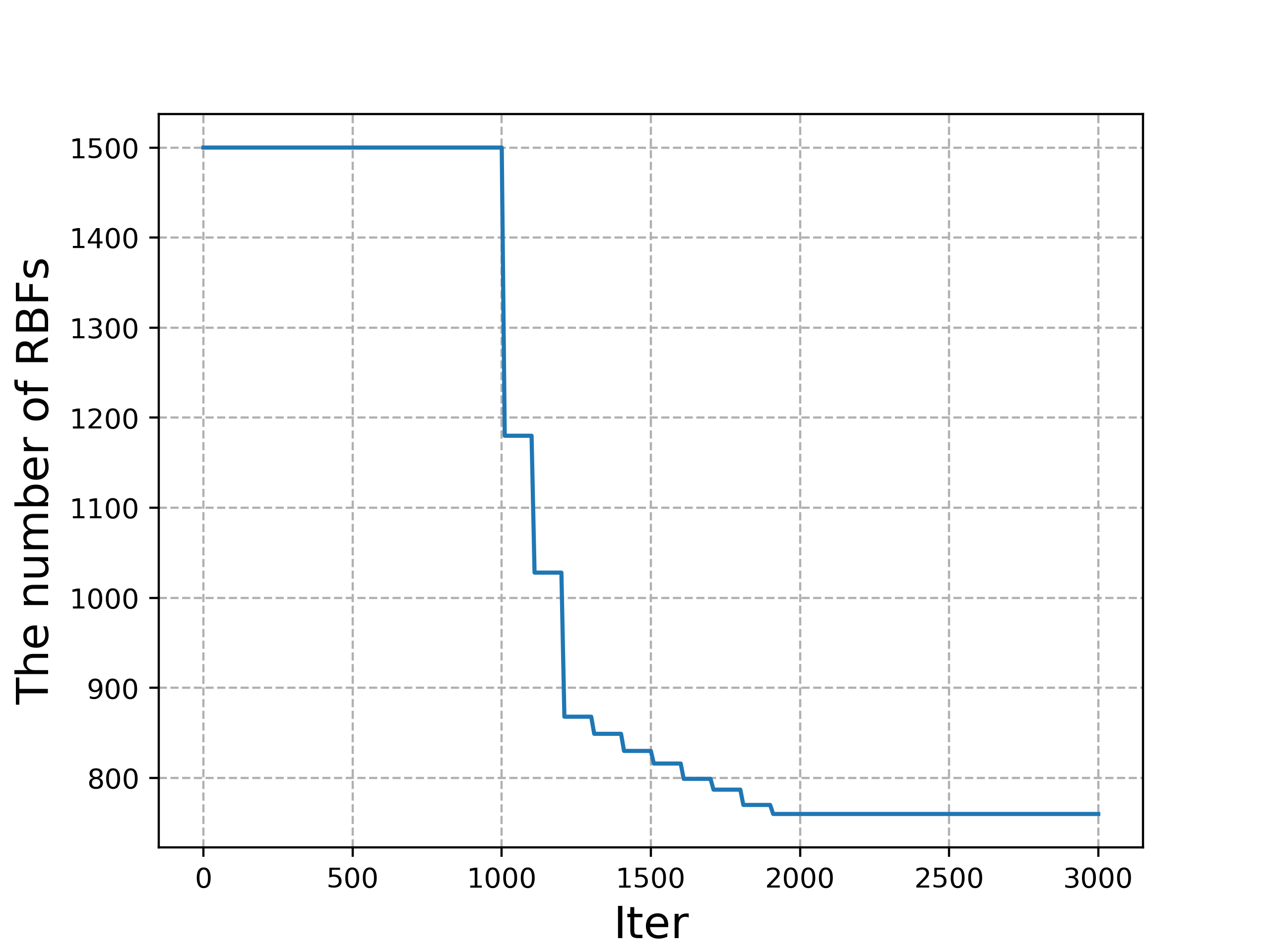}
	}
	\caption{Three relative errors and the number of RBFs in terms of the iteration number when $\varepsilon$=0.002 in Example 3.}
	\label{figexam3}
\end{figure}

\paragraph{Example 3}
Consider another two-scale function
\begin{eqnarray}
	a^{\varepsilon}(x)=2+\sin(2\pi x+\frac{2\pi x}{\varepsilon}).
	\label{exam3}
\end{eqnarray}
Fig. \ref{figexam3} records the training processes of three relative errors and the number of RBFs when $\varepsilon=0.002$. According to Fig. \ref{figexam3}, we can observe that three relative errors stabilize around the 500th iteration, but they are still large at this point. When
the regularization term was added in the loss function at around the 1000th iteration, the three relative errors all decreased
significantly, which means that the $ell_1$ regularization term can improve the generalization performance of the network.
Table \ref{tabexam3} records the results of three relative errors and the number of RBFs in the final solution for Example 3. It is recognized that our method works well for two-scale coefficients. The accuracy and the number of RBFs in the final solution are also similar to the first two examples.
\begin{table}[H]
	\centering
	\begin{tabular}{ccccc}
		\hline
		$\varepsilon$ & $N$ & $err_2$ & $err_\infty$ & $err_{H_1}$  \\ \hline	
		\rowcolor{gray!40}
		0.5  &18  & 4.850e-5  & 1.336e-4 & 2.319e-4 \\
		0.1  &40  & 6.539e-5  & 3.710e-4 & 5.730e-4 \\
		\rowcolor{gray!40}
		0.05 & 69 & 6.778e-5  & 4.674e-4 & 4.249e-4 \\
		0.01 &178 & 1.086e-4  & 3.675e-4 & 1.414e-3 \\
		\rowcolor{gray!40}
		0.005& 379& 1.431e-4  & 2.798e-4 & 2.134e-3 \\
		0.002& 760& 9.549e-5  & 3.698e-4 & 1.500e-3 \\\hline
	\end{tabular}
	\caption{Relative $L_2$, $L_{\infty}$ and $H_1$ errors and the number of basis functions in the final solution for Example 3. }
	\label{tabexam3}
\end{table}
\paragraph{Example 4}
Consider a discontinuous and periodic function
\begin{equation}
	\begin{aligned}
		a^{\varepsilon}(x)=\left\{\begin{aligned}
			1 \quad x\in [0,\frac{\varepsilon}{2}) \\
			10 \quad x\in [\frac{\varepsilon}{2},\varepsilon)
		\end{aligned}
		\right..
	\end{aligned}
	\label{exam4}
\end{equation}
Fig. \ref{figexam4_1} plots the solution and its derivative obtained by SRBFNN and FDM when $\varepsilon=0.05$, 0.002, respectively. It is observed that SRBFNN has a good approximation for both the solution and the derivative even when the coefficient is discontinuous.
\begin{small}
	\begin{figure}[H]
		\centering
		\subfigure[$\varepsilon=0.05$]{ 
			\includegraphics[align=c,width=0.4\columnwidth]{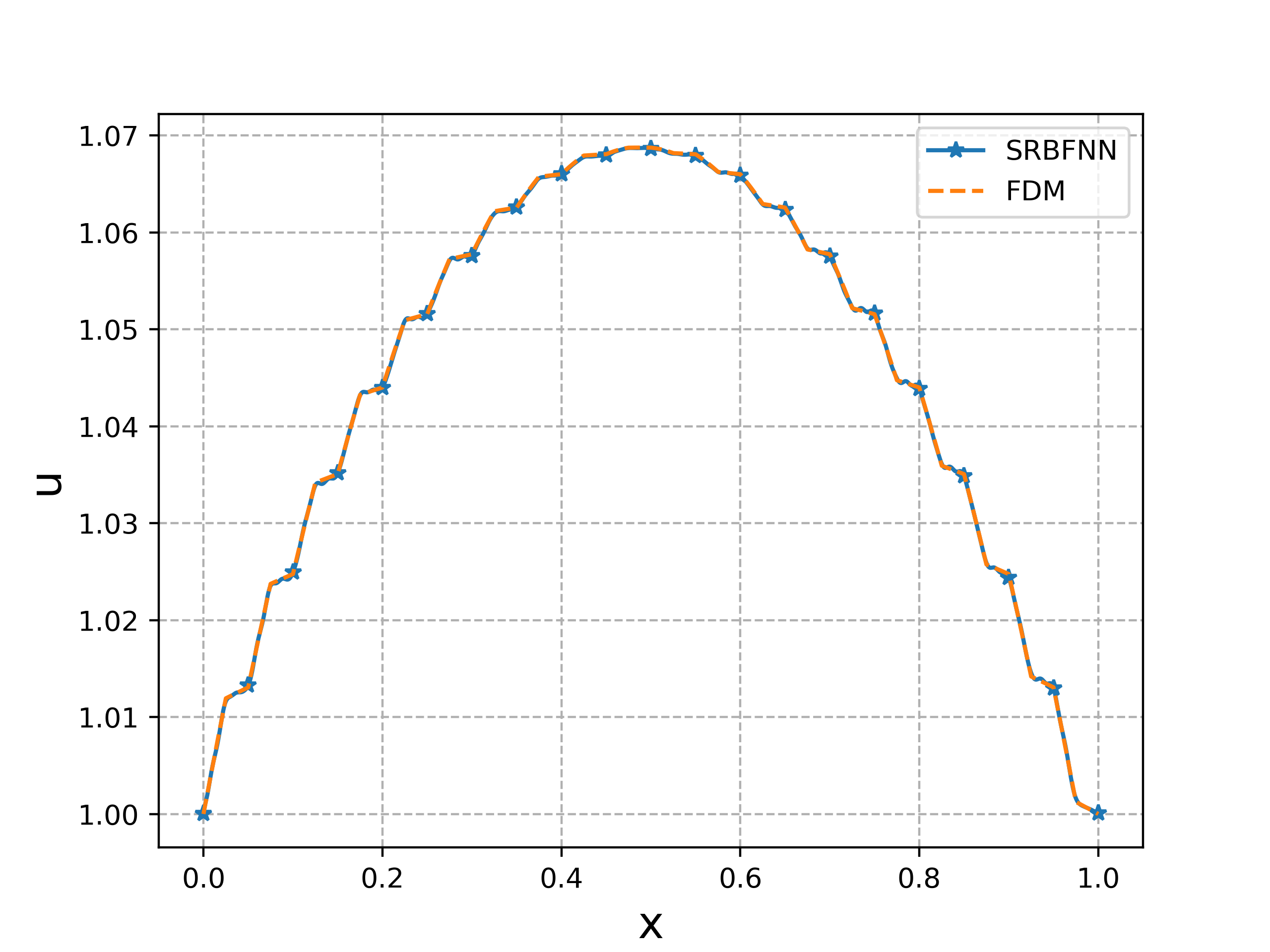}
		}
		\hfill 
		\subfigure[$\varepsilon=0.05$]{ 
			\includegraphics[align=c,width=0.4\columnwidth]{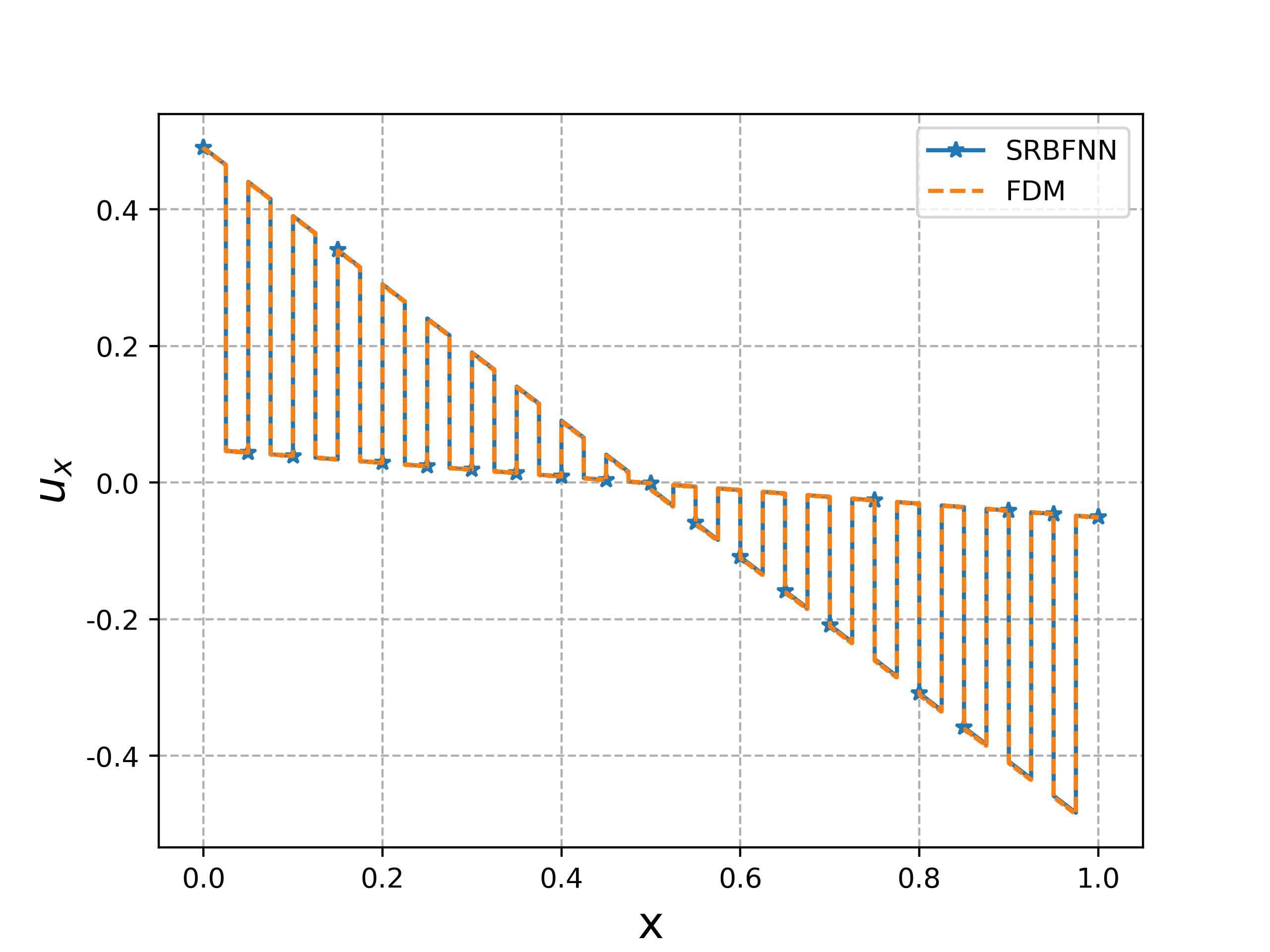}
		}
		\subfigure[$\varepsilon=0.002$]{ 
			\includegraphics[align=c,width=0.4\columnwidth]{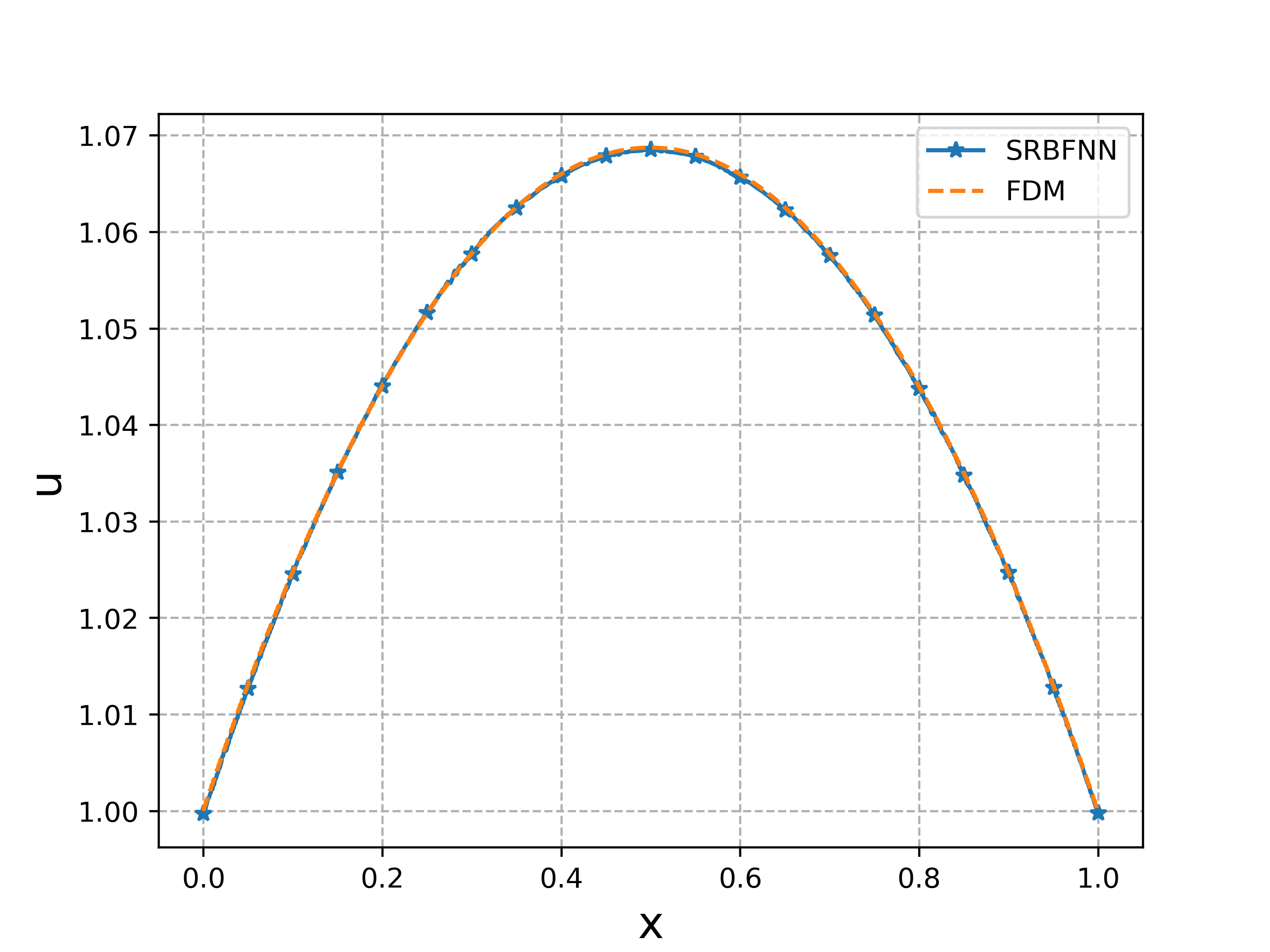}
		}
		\hfill 
		\subfigure[$\varepsilon=0.002$]{ 
			\includegraphics[align=c,width=0.4\columnwidth]{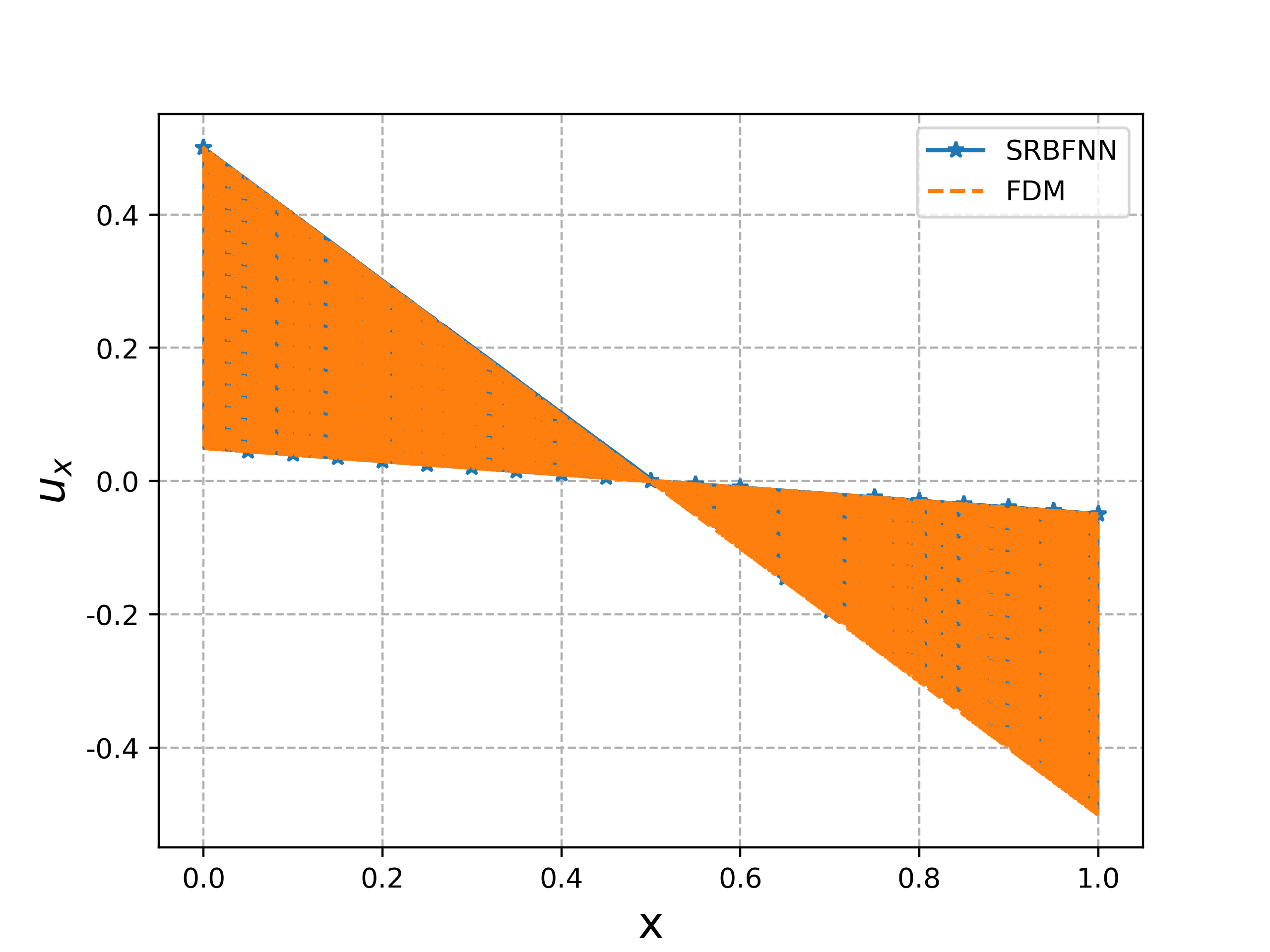}
		}
		\caption{Numerical solution and its derivative obtained by SRBFNN and FDM when $\varepsilon$=0.05, 0.002 for Example 4. }
		\label{figexam4_1}
	\end{figure}
\end{small}
The relative $L_2$ and $H_1$ errors for five methods are shown in Table \ref{tabexam4}. The network architectures for MscaleDNN, DRM, DGM and PINN have the same settings as those used in Example 2. 
\begin{table}[H]
	\centering
	\resizebox{.9\textwidth}{!}{
		\begin{tabular}{c|ccccc|ccccc}
			\toprule
			& &\multicolumn{3}{c}{$err_2$}& & &\multicolumn{3}{c}{$err_{H_1}$}&\\ 
			
			$\varepsilon$ &SRBFNN&MscaleDNN&DRM&DGM&PINN  &  SRBFNN&MscaleDNN&DRM&DGM&PINN\\ \hline	
			\rowcolor{gray!40}
			0.5&3.673e-4&1.241e-3&1.668e-2&1.544e-2&1.698e-2 & 6.391e-3&2.353e-2&9.213e-2&1.121e-1&1.152e-1\\
			0.1&1.342e-4&1.375e-2&3.060e-2&3.712e-2&3.871e-2 & 5.174e-3&8.739e-2&1.583e-1&1.739e-1&1.740e-1\\
			\rowcolor{gray!40}
			0.05&1.185e-4&2.396e-2&3.092e-2&3.771e-2&3.898e-2& 7.461e-3&1.273e-1&1.597e-1&1.744e-1&1.746e-1\\
			0.01&3.359e-4&3.127e-2&3.108e-2&3.622e-2&3.811e-2& 1.765e-2&1.579e-1&1.595e-1&1.735e-1&1.737e-1\\
			\rowcolor{gray!40}
			0.005&6.714e-4&3.189e-2&3.104e-2&3.839e-2&3.846e-2& 2.458e-2&1.587e-1&1.585e-1&1.730e-1&1.731e-1\\
			0.002&2.218e-4&3.201e-2&3.107e-2&3.661e-2&3.794e-2& 3.832e-2&1.560e-1&1.557e-1&1.698e-1&1.705e-1\\ 
			\bottomrule
		\end{tabular}
	}
	\caption{Relative $L_2$ and $H_1$ errors of SRBFNN, MscaleDNN, DRM, DGM, and PINN for Example 4. The penalty coefficients of the boundary condition and the maximum number of iterations are the same as Example 2 for MscaleDNN, DRM, DGM, and PINN.}
	\label{tabexam4}
\end{table}

As depicted in Table \ref{tabexam4}, SRBFNN obtains more accurate results than MscaleDNN, DRM, DGM, and PINN. Again, in the same settings, all the errors in Example 4 are larger than those in Example 2 and Example 3, while the discontinuity has less impact on SRBFNN. These results indicate that SRBFNN can handle multiscale problems with  discontinuous coefficients.

\begin{figure}[H]
	\centering
	\includegraphics[align=c,width=0.45\columnwidth]{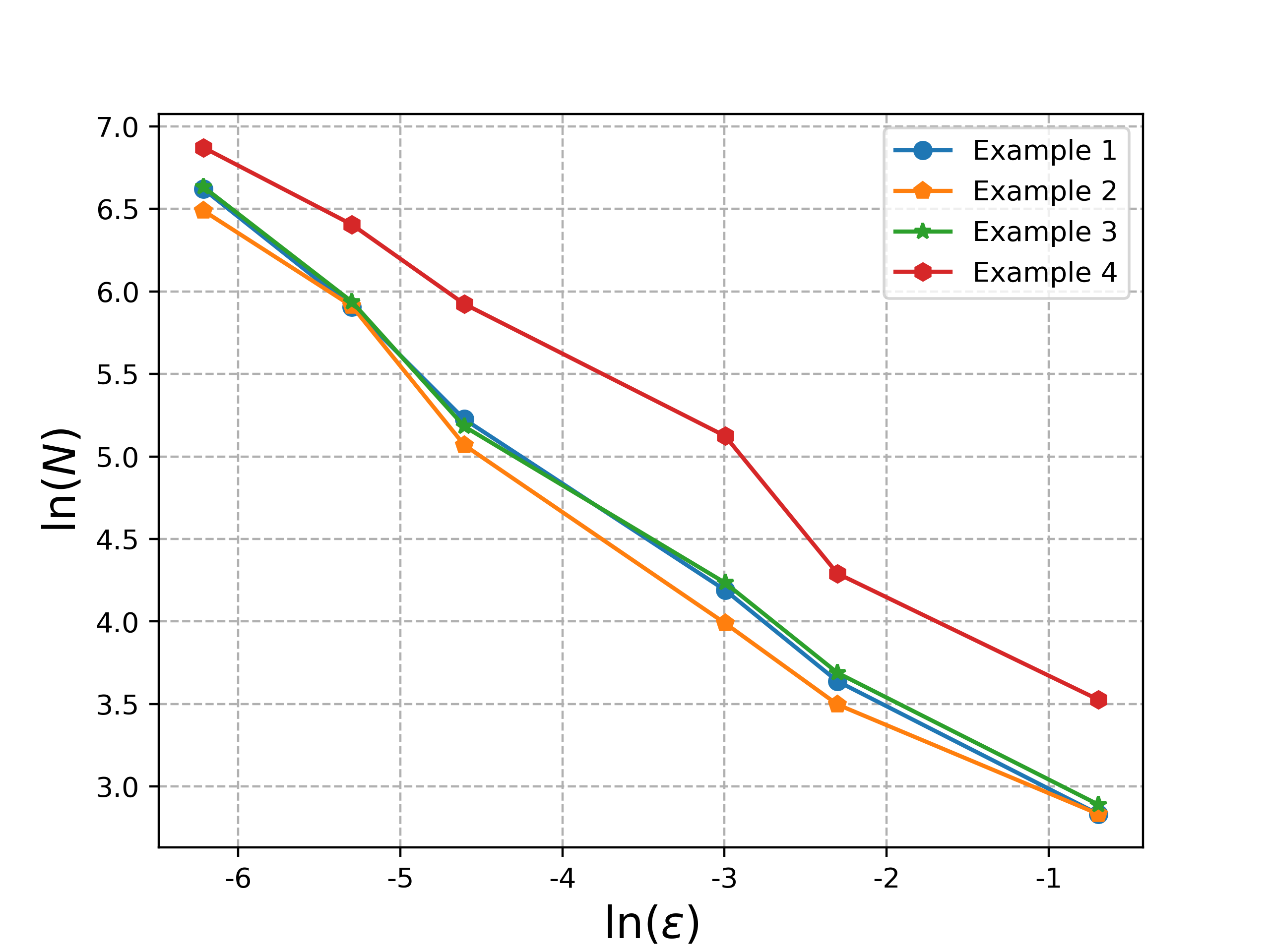}
	\caption{The number of basis functions needed in the final solution as a function of $\varepsilon$ for Examples 1-4.}
	\label{fig1d}
\end{figure}

Next, we plot the number of RBFs needed in the final solution and $\varepsilon$ for Examples 1-4 in Fig. \ref{fig1d}. It is observed that ln$(N)$ and ln$(\varepsilon^{-1})$ is approximatively linear and more RBFs are needed in Example 4 than in other examples. Table \ref{tab7} estimates the slope by using the least squares method. From Table \ref{tab7}, we see that ln$(N)\sim$ 0.67ln$(\varepsilon^{-1}) $, which implies  $N=\mathcal{O}(\varepsilon^{0.67})$.
\begin{table}[H]
	\centering
	\begin{tabular}{cccc}
		\hline
		Example 1 & Example 2 & Example 3 & Example 4 \\ \hline
		\rowcolor{gray!40}
		0.693     & 0.687     & 0.683     & 0.621     \\ \hline
	\end{tabular}
	\caption{The least squares estimation of slopes of ln$(N)$ with respect to ln$(\varepsilon^{-1})$ for Examples 1-4. }
	\label{tab7}
\end{table}

\subsection{Numerical experiments in 2D }
\textbf{Setting:}\label{2dset}
In all two-dimensional cases, the initial learning rate is set to $0.1$ and reduces by 1/10 every 40 iterations until it is less than $10^{-5}$. The training data is equidistantly sampled with mesh size $h$=0.002. The mesh size of FDM is set to 0.0005. The $batchsize$ is 1024 in the domain and $512\times4$ on the boundary.
$MaxNiter=300$ , $SparseNiter=250$, $tol_1=0.1$, $tol_2=10^{-5}$, $Check_{iter}=10$. The hyperparameters $\lambda_{2}$ and $\lambda_{3}$ are initialized to 20.0, 0.001 respectively. The solution is more oscillatory as $\varepsilon$ decreases, and this makes the magnitude of derivatives larger. So to balance the effect of different loss terms, the penalty coefficient $\lambda_{1}$ is gradually reduced when $\varepsilon$ goes to smaller as follows.
\begin{equation}
	\begin{aligned}
		\lambda_{1}=\left\{\begin{aligned}
			0.1 &\mbox{, if } \varepsilon=0.5,0.2,0.1 \\
			0.02 &\mbox{, if } \varepsilon=0.05,0.02\\
			0.004 &\mbox{, if } \varepsilon=0.01
		\end{aligned}
		\right..
	\end{aligned}
\end{equation}

Consider the problems in the domain $\Omega=[0,1]^2$ with six scales, $\varepsilon=$0.5, 0.2, 0.1, 0.05, 0.02, 0.01 for all 2D examples.
The initial number of RBFs is 1000, 1000, 2000, 5000, 15000, and 30000, respectively.

\begin{figure}[H]
	\centering
	\subfigure[Relative errors]{ 
		\includegraphics[align=c,width=0.45\columnwidth]{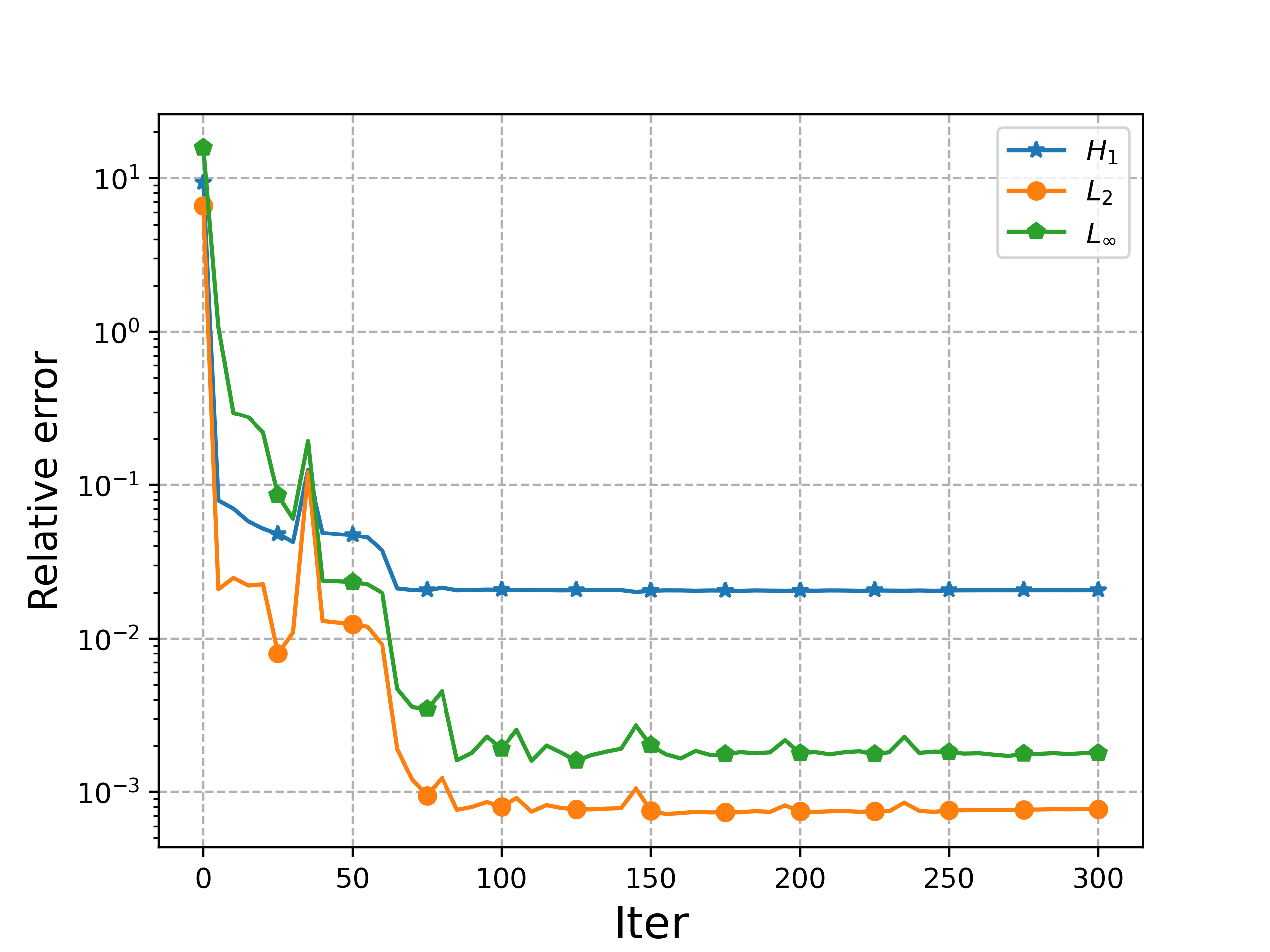}
	}
	\hfill %   ???   ?   
	\subfigure[The number of RBFs]{ 
		\includegraphics[align=c,width=0.45\columnwidth]{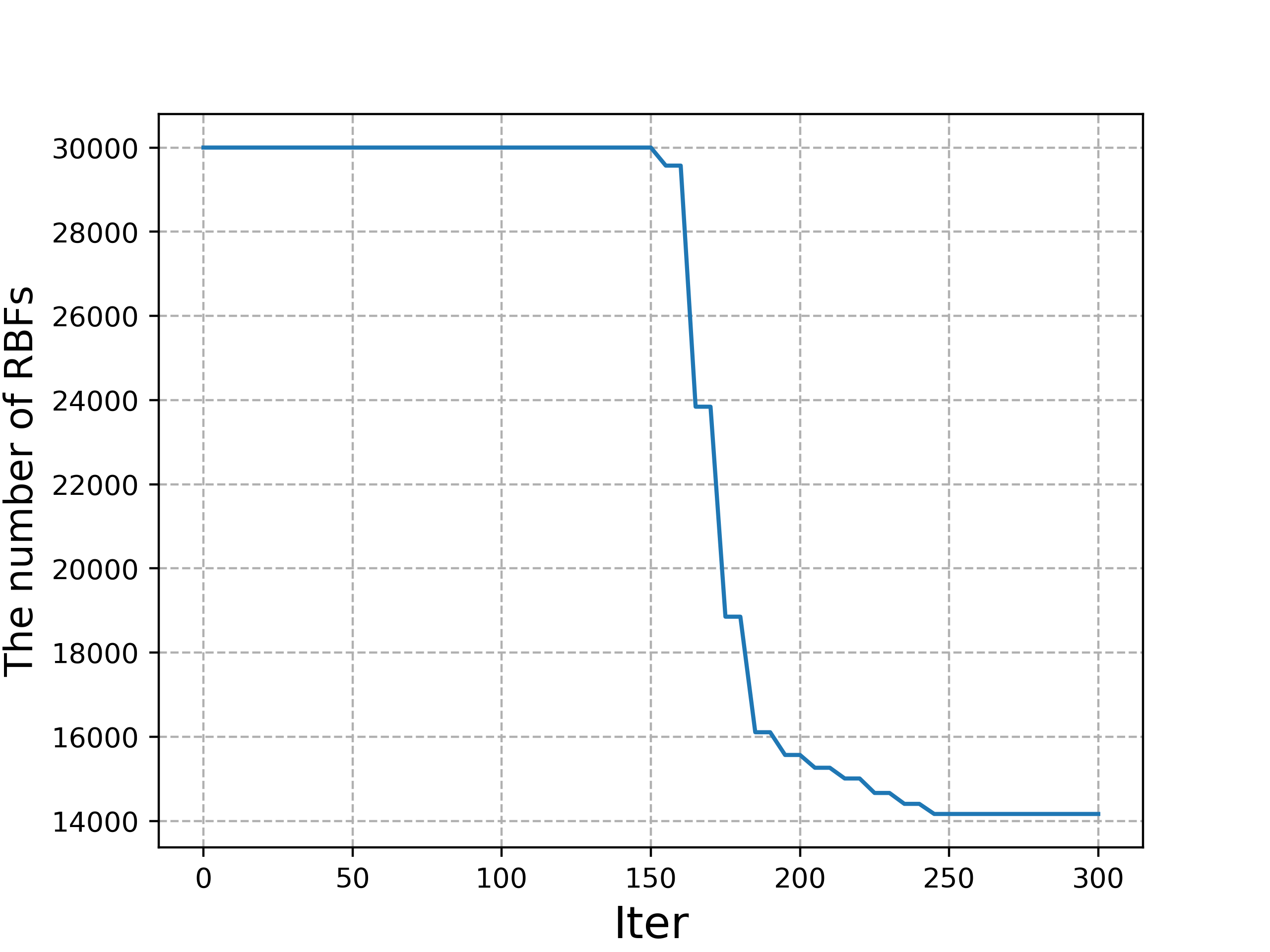}
	}
	\caption{Training process of three relative errors and the number of RBFs when $\varepsilon$= 0.01 for Example 5. }
	\label{figexam5_1}
\end{figure}

\paragraph{Example 5}
$a^{\varepsilon}(x,y)$ is defined by
\begin{equation}
	a^{\varepsilon}(x,y)=2+\sin(\frac{2\pi (x+y)}{\varepsilon} ),
	\label{eqexam5}
\end{equation}
where $f(x,y)=-1.0$, $g(x,y)=1.0$. 
Fig. \ref{figexam5_1} plots the training processes of three relative errors and the number of RBFs when $\varepsilon=0.01$ for Example 5. One can see that SRBFNN achieves a good accuracy around the $100th$ iteration and the number of basis functions no longer decreases at the $250th$ iteration.

\begin{figure}[H]
	\centering
	\subfigure[$\varepsilon=0.5$]{ 
		\includegraphics[align=c,width=0.30\columnwidth]{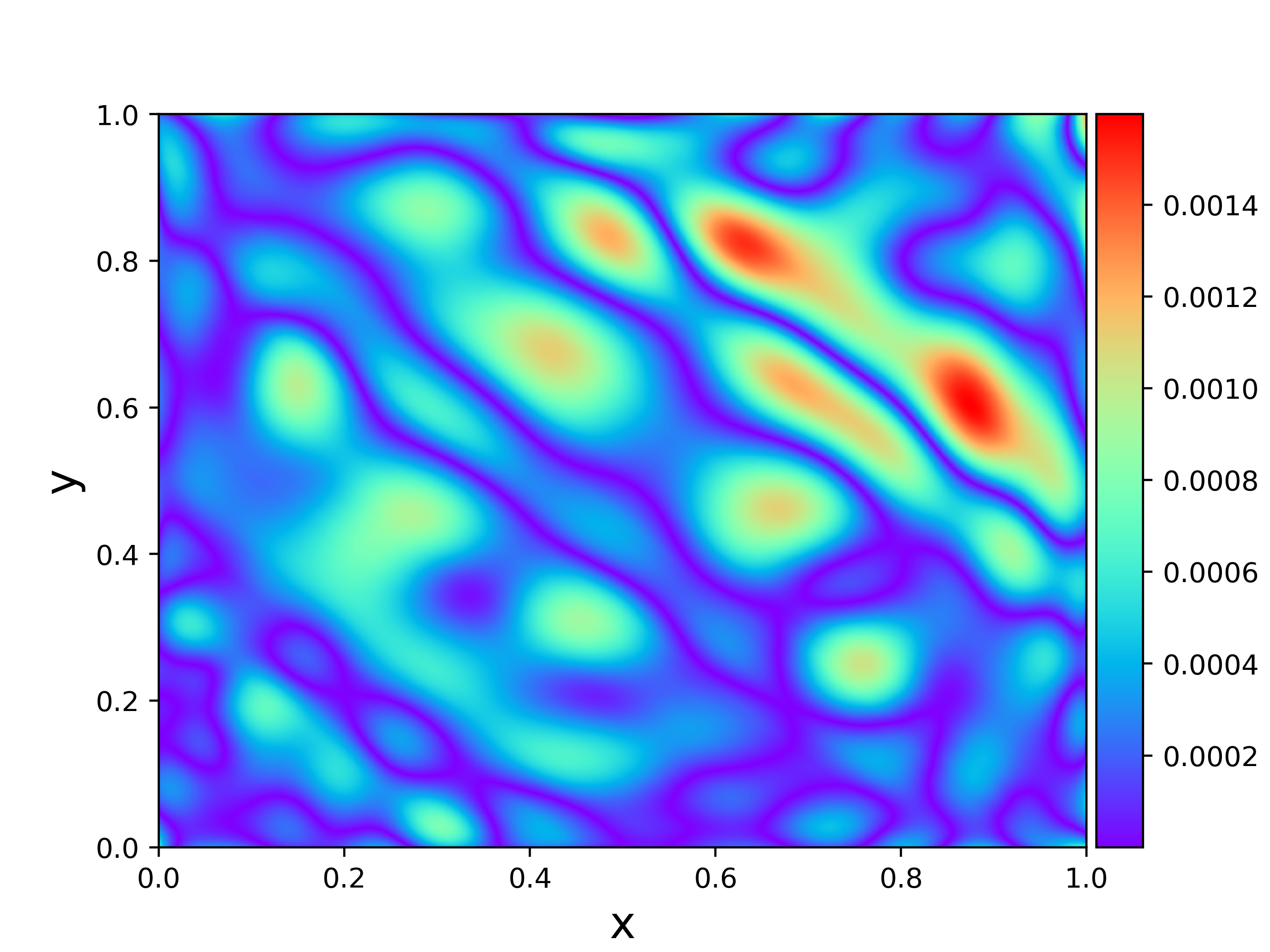}
	}
	\hfill %   ???   ?   
	\subfigure[$\varepsilon=0.2$]{ 
		\includegraphics[align=c,width=0.30\columnwidth]{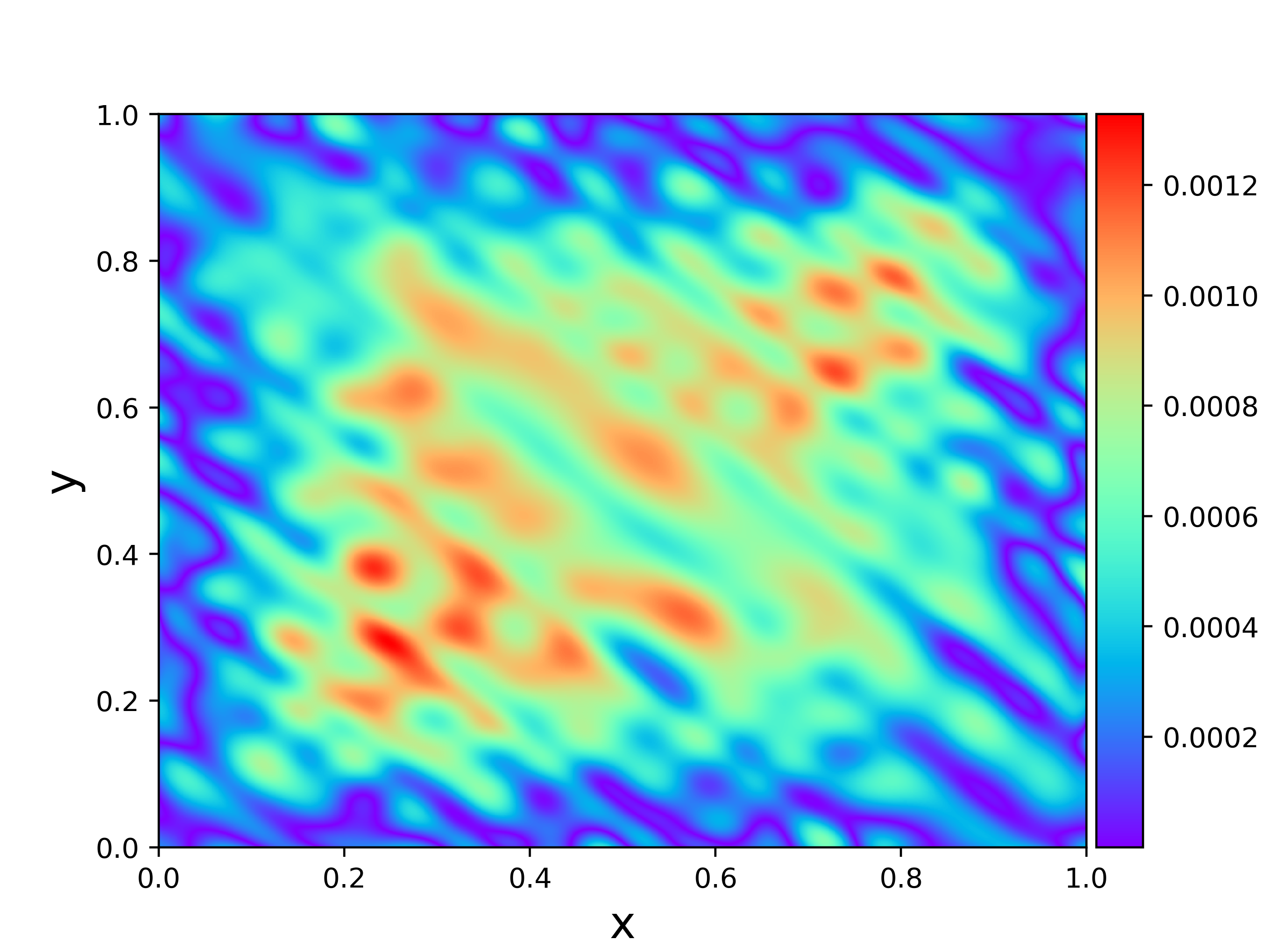}
	}
	\hfill %   ???   ?   
	\subfigure[$\varepsilon=0.1$]{ 
		\includegraphics[align=c,width=0.30\columnwidth]{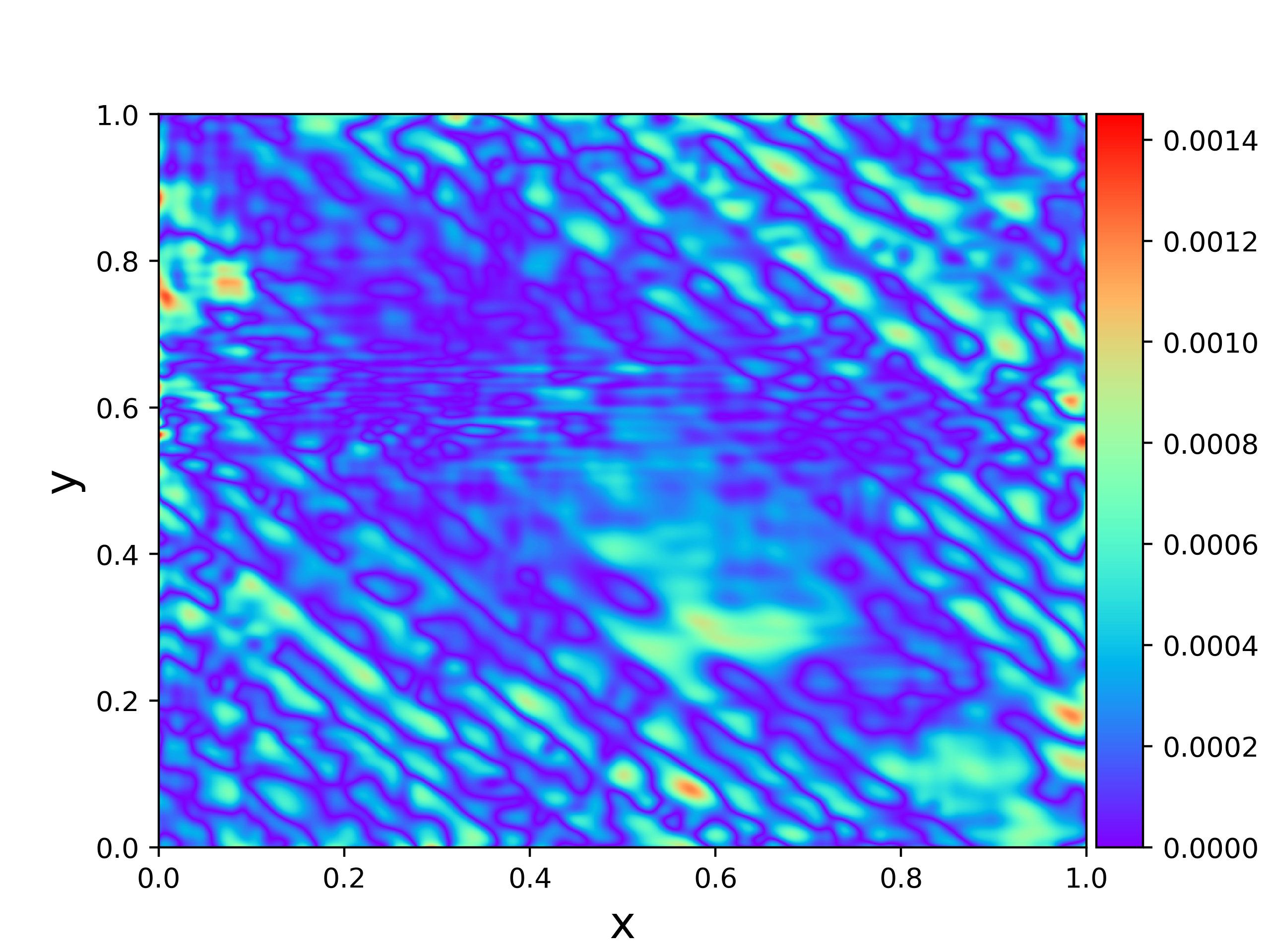}
	}
	\hfill %   ???   ?   
	\subfigure[$\varepsilon=0.05$]{ 
		\includegraphics[align=c,width=0.30\columnwidth]{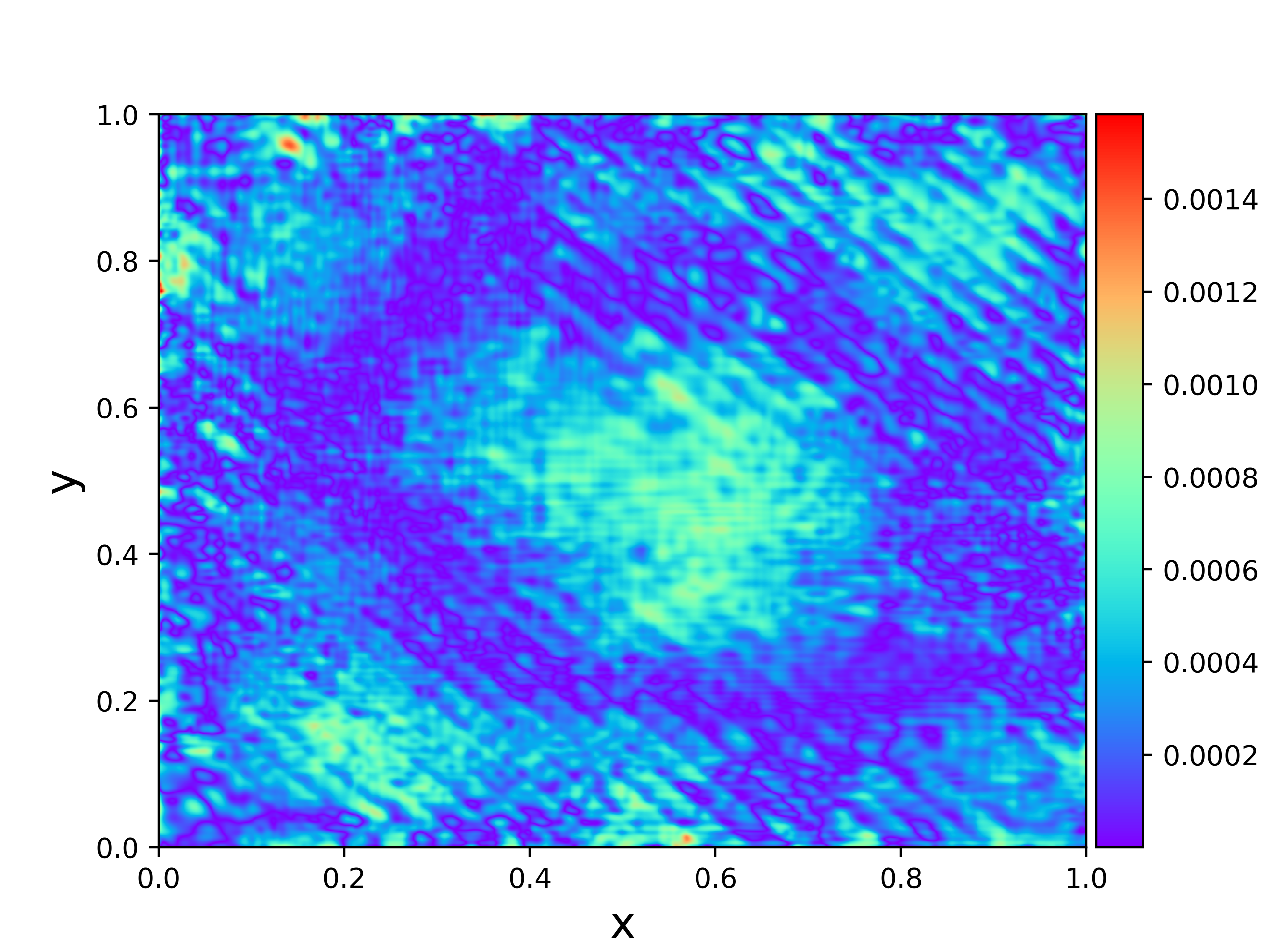}
	}
	\hfill %   ???   ?   
	\subfigure[$\varepsilon=0.02$]{ 
		\includegraphics[align=c,width=0.30\columnwidth]{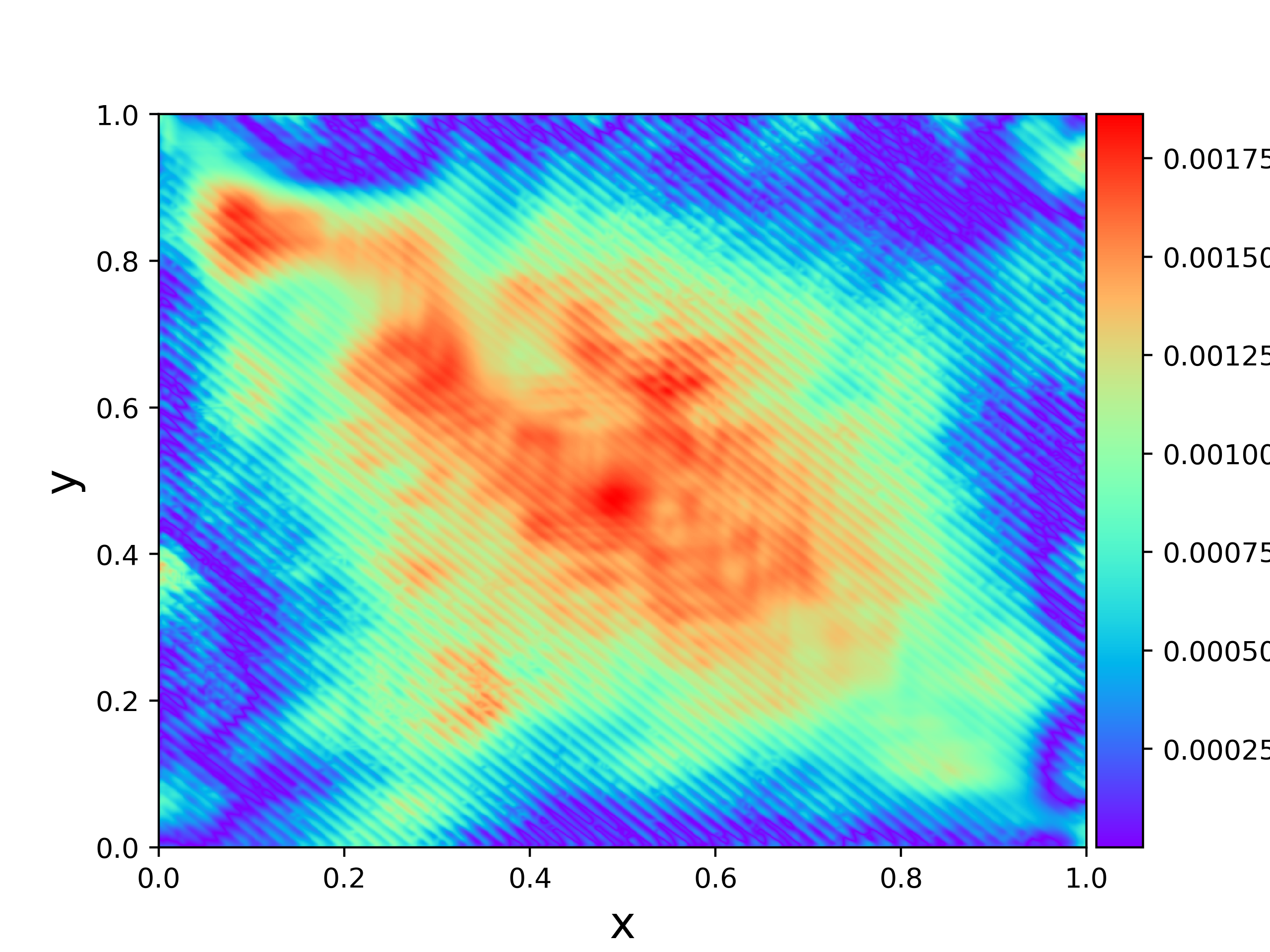}
	}
	\hfill %   ???   ?   
	\subfigure[$\varepsilon=0.01$]{ 
		\includegraphics[align=c,width=0.30\columnwidth]{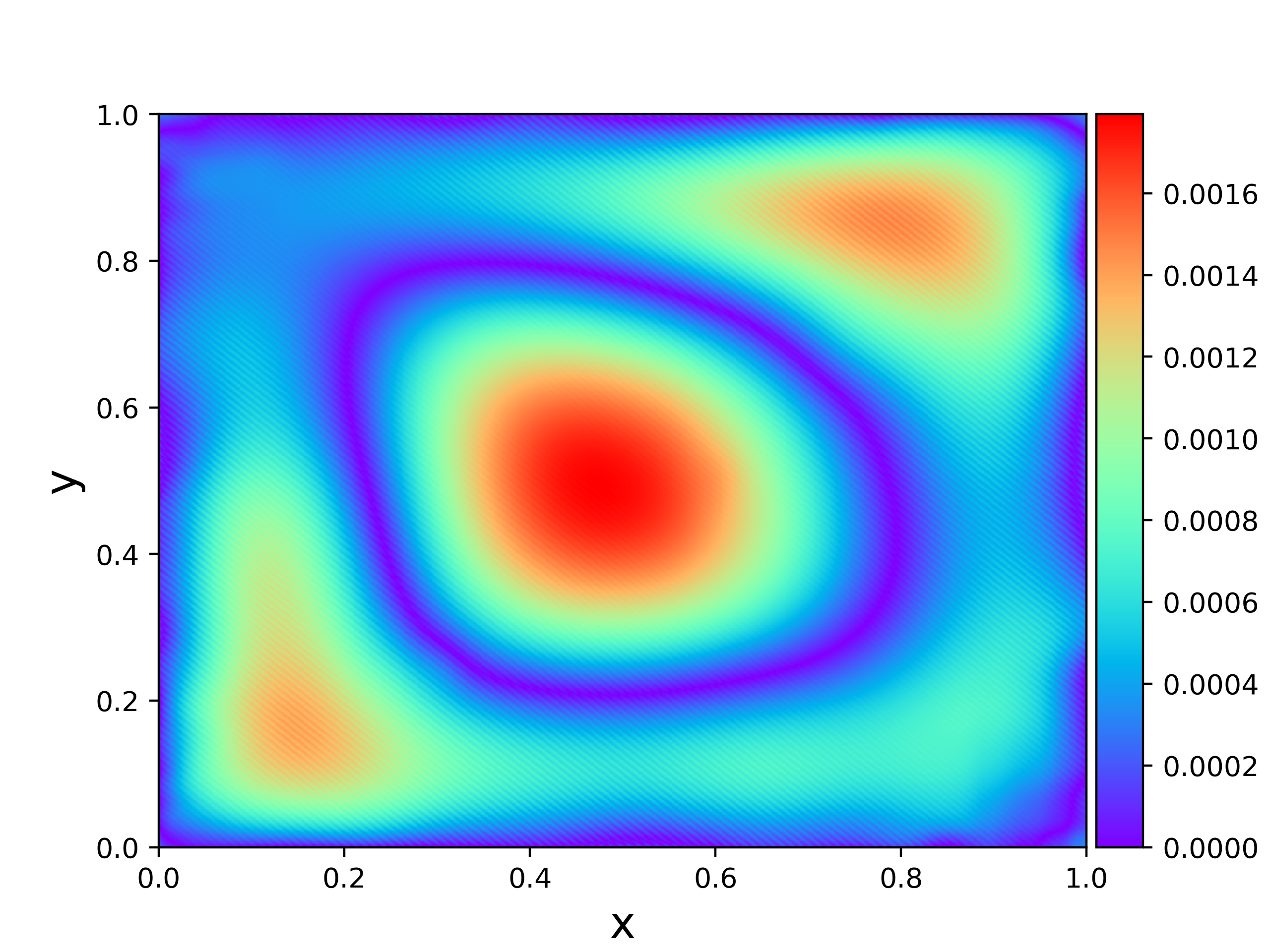}
	}
	\hfill %   ???   ?   
	\caption{The absolute point-wise errors for different $\varepsilon$ for Example 5.}
	\label{figexam5}
\end{figure}

A detailed assessment of the predicted solution is presented in Fig. \ref{figexam5}. In particular, it shows the absolute point-wise errors $|u^S-u^F|$ at different scales. We can observe from Fig. \ref{figexam5} that the absolute point-wise errors in most areas are the order of O($10^{-3}$), which means our method provides a good approximation for the solution. Then three relative errors and the number of basis functions in the final solution are recorded in Table \ref{tabexam5}. One can see that the number of basis functions grows much faster than that in 1D when $\varepsilon$ becomes smaller. In addition, the accuracy of relative $L_2$ and $L_{\infty}$ errors is over one order of magnitude than that in 1D.
\begin{table}[H]
	\centering
	\begin{tabular}{ccccc}
		\hline
		$\varepsilon$ & $N$ & $err_2$ & $err_\infty$ & $err_{H_1}$  \\ \hline	
		\rowcolor{gray!40}
		0.5    & 42  & 4.798e-4  & 1.596e-3 & 5.442e-3 \\
		0.2    & 146 & 6.080e-4  & 1.328e-3 & 1.721e-2 \\
		\rowcolor{gray!40}
		0.1    & 531 & 3.164e-4  & 1.453e-3 & 1.772e-2 \\
		0.05   & 2121& 3.561e-4  & 1.582e-3 & 2.684e-2 \\
		\rowcolor{gray!40}
		0.02   & 6272& 9.640e-4  & 1.862e-3 & 3.189e-2 \\
		0.01   &14164& 7.719e-4  & 1.793e-3 & 2.070e-2 \\\hline
	\end{tabular}
	\caption{Relative $L_2$, $L_{\infty}$ and  $H_1$ errors and the number of basis functions in the final solution for Example 5. }
	\label{tabexam5}
\end{table}

\begin{figure}[H]
	\centering
	\includegraphics[align=c,width=0.45\columnwidth]{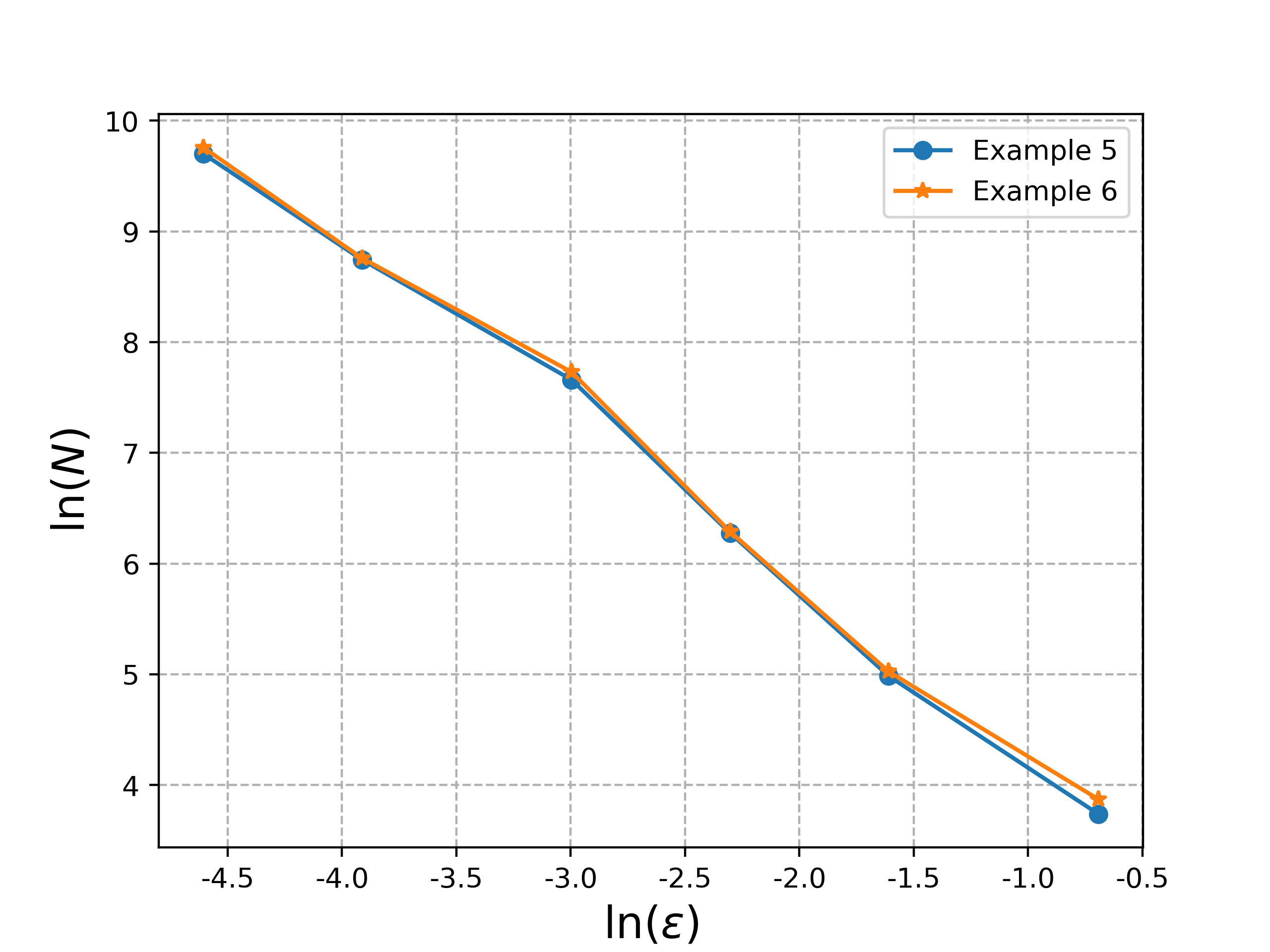}
	\caption{The number of basis functions in terms of $\varepsilon$ in the log-log scale for Example 5 and 6.}
	\label{figexam6_1}
\end{figure}

\paragraph{Example 6 \cite{ming2006}}
Consider $a^\varepsilon(x,y)$ is a function with two scales
\begin{eqnarray}
	a^{\varepsilon}(x,y)=\frac{1.5+\sin(\frac{2\pi x}{\varepsilon})}{1.5+\sin(\frac{2\pi y}{\varepsilon} )}+\frac{1.5+\sin(\frac{2\pi y}{\varepsilon})}{1.5+\cos(\frac{2\pi x}{\varepsilon})}+\sin(4x^2y^ 2)+1,
	\label{exam6}
\end{eqnarray}
where $f(x,y)=-10.0$, $g(x,y)=0.0$. We conduct a comparison for SRBFNN with MscaleDNN, DRM, DGM, and PINN. The DRM and DGM contain four residual blocks. The PINN has nine linear layers. The hidden layers are 64, 64, 128, 128, 256, 256 when $\varepsilon$=0.5, 0.2, 0.1, 0.05, 0.02, 0.01 for DRM, DGM and PINN. The hidden layers of MscaleDNN are (1000,400,300,300,200,100,100) for all $\varepsilon$.  $batchsize$ is 1024 in the domain and $512\times4$ on the boundary for these four methods. The initial learning rate is 0.001 and reduces by 1/10 every 100 iterations until it is less than $10^{-5}$. We set the maximum number of iterations as 1000 for MsacleDNN, DRM, DGM, and PINN. Table \ref{tabexam6} shows the relative $L_2$ and $H_1$ errors in terms of different scales for five methods. 
\begin{table}[H]
	\centering
	\resizebox{.9\textwidth}{!}{  % Here 1/2
		\begin{tabular}{c|ccccc|ccccc}
			\toprule
			& &\multicolumn{3}{c}{$err_2$}& &&\multicolumn{3}{c}{$err_{H_1}$}  &\\ 
			$\varepsilon$ &SRBFNN&MscaleDNN&DRM&DGM&PINN & SRBFNN&MscaleDNN&DRM&DGM&PINN\\ \hline	
			\rowcolor{gray!40}
			0.5 &1.319e-3&7.885e-3&2.062e-2&6.615e-3&8.812e-3 & 1.276e-2&5.979e-2&1.074e-1&1.005e-1&1.048e-1\\
			0.2 &2.116e-3&8.467e-3&2.280e-2&7.393e-2&1.658e-1 & 1.650e-2&9.922e-2&1.654e-1&1.182e-1&0.303e-1\\
			\rowcolor{gray!40}
			0.1 &3.049e-3&8.229e-3&2.601e-2&7.531e-1&4.333e-1 & 3.497e-2&1.116e-1&2.071e-1&7.626e-1&7.164e-1\\
			0.05&3.150e-3&3.096e-2&2.681e-2&9.676e-1&9.133e-1 & 7.774e-2&2.164e-1&2.182e-1&9.466e-1&9.231e-1\\
			\rowcolor{gray!40}
			0.02&3.395e-3&2.066e-2&3.546e-2&9.965e-1&1.000e-0 & 9.972e-2&2.019e-1&2.210e-1&9.958e-1&9.964e-1\\
			0.01&2.481e-3&2.589e-2&3.667e-2&9.982e-1&1.009e-0 & 1.125e-1&2.078e-1&2.392e-1&9.993e-1&9.997e-1\\ 
			\bottomrule
		\end{tabular}
	}
	\caption{Relative $L_2$ and $H_1$ errors of SRBFNN, MscaleDNN, DRM, DGM, and PINN for Example 6. The penalty coefficient of the boundary condition for MscaleDNN and DRM is 500.0, and is 10.0 for PINN and DGM. }
	\label{tabexam6}
\end{table}
As depicted in Table \ref{tabexam6}, SRBFNN reaches relatively lower errors than MscaleDNN, DRM, DGM, and PINN. PINN and DGM perform well when $\varepsilon$ is big but the errors are relatively larger when $\varepsilon$ is small. Additionally, one can observe that MscaleDNN andd DRM outperform DGM and PINN, while MscaleDNN, DGM, PINN outperform DRM in one-dimensional cases. This observation is similar to the conclusion in \cite{chen2020comparison}. Then we compare the number of parameters and the average time per iteration in all methods. 
\begin{table}[H]
	\centering
	\resizebox{.9\textwidth}{!}{  % Here 1/2
		\begin{tabular}{c|ccccc|ccccc}
			\toprule
			&\multicolumn{5}{c}{Parameters}&\multicolumn{5}{c}{Average time per iteration (Seconds) }\\ 
			$\varepsilon$ &SRBFNN&MscaleDNN& DRM & DGM   & PINN&SRBFNN&MscaleDNN& DRM & DGM   & PINN\\ \hline	
			\rowcolor{gray!40}
			0.5  &5240&704501 &33537 &33537 &29377&11.102 &9.731  &6.948  &11.725  & 10.087 \\
			0.2  &5760&704501 &33537 &33537 &29377&10.747 &10.315  &6.722  &11.796  &10.111  \\
			\rowcolor{gray!40}
			0.1  &12690&704501&132609&132609&116097 &10.889 &9.885  &6.722  &11.821  &10.233  \\
			0.05 &36375&704501&132609&132609&116097 &15.310 &9.713  &6.751  &11.783  &10.185  \\
			\rowcolor{gray!40}
			0.02 &106810&704501&527361&527361&461569&37.869 &9.821 &6.827  &11.704  &10.032  \\
			0.01 &236425&704501&527361&527361&461569&74.965 &8.231  &6.756  &11.942  &10.080  \\ \bottomrule
		\end{tabular}
	}
	\caption{The network parameters and the average running time per iteration in different methods for Example 6.}
	\label{tabexam6_1}
\end{table}
From Table \ref{tabexam6_1},  it
is recognized that SRBFNN achieves better approximation than MscaleDNN, DRM, DGM, and PINN with fewer parameters. However, the average time per iteration for SRBFNN is larger than MscaleDNN, DRM, DGM, and PINN when $\varepsilon$ is small, but SRBFNN requires fewer iterations to converge. Fig. \ref{figexam6_1} plots the final number of RBFs in terms of $\varepsilon$. Analogously, we use the least squares method to estimate the slope of ln($N$) as a function of ln($\varepsilon$) for Example 5 and Example 6 in Table \ref{tabexam6_2}. One can observe from Fig. \ref{figexam6_1} and Table \ref{tabexam6_2} that ln($N)\sim 1.55$ln($\varepsilon^{-1}$), which implies $N=\mathcal{O}(\varepsilon^{-0.78* 2})$. This observation, together with the results in aforementioned 1D examples can be written in a more general form: $N=\mathcal{O}(\varepsilon^{-\tau n})$, $N$ is the number of RBFs in final solution, $n$=1, 2 is the dimensionality and $\tau$ is smaller than $1$.
\begin{table}[H]
	\centering
	\begin{tabular}{ll}
		\hline
		Example 5 &Example 6 \\ \hline
		\rowcolor{gray!40}
		1.560	  & 1.544      \\ \hline
	\end{tabular}
	\caption{The least squares estimation of slope about ln$(N)$ and ln$(\frac{1}{\varepsilon})$ for Example 5 and Example 6. }
	\label{tabexam6_2}
\end{table}

\paragraph{Example 7 (\cite{ming2006})}
Consider $a^\varepsilon(x,y)$ having six scales
\begin{equation}
	\begin{split}
		a^{\varepsilon}(x,y)=\frac{1}{6}&(\frac{1.1+\sin(2\pi x/\varepsilon_1)}{1.1+\sin(2\pi y/\varepsilon_1 )}+\frac{1.1+\sin(2\pi y/\varepsilon_2)}{1.1+\cos(2\pi x/\varepsilon_2)}+\frac{1.1+\cos(2\pi x/\varepsilon_3)} {1.1+\sin(2\pi y/\varepsilon_3)}+\nonumber\\
		&\frac{1.1+\sin(2\pi y/\varepsilon_4)}{1.1+\cos(2\pi x/\varepsilon_4)}+\frac{1.1+\cos(2\pi x/\varepsilon_5)}{1.1 +\sin(2\pi y/\varepsilon_5)}+\sin(4x^2y^2)+1 ),
	\end{split}
	\label{exam7}
\end{equation}
where $f(x,y)=-10.0$, $g(x,y)=0.0$, $\varepsilon_1=1/5,\varepsilon_2=1/13,\varepsilon_3=1/17,\varepsilon_4=1/31,\varepsilon_5=1/65$. Particularly, the initial learning rate is set to $0.1$ and reduces to 1/10 every 40 iterations until it is less than $10^{-6}$. The penalty parameter $\lambda_{1}$ is 0.02 and the initial number of RBFs is 30000.

\begin{figure}[H]
	\centering
	\includegraphics[width=0.4\linewidth]{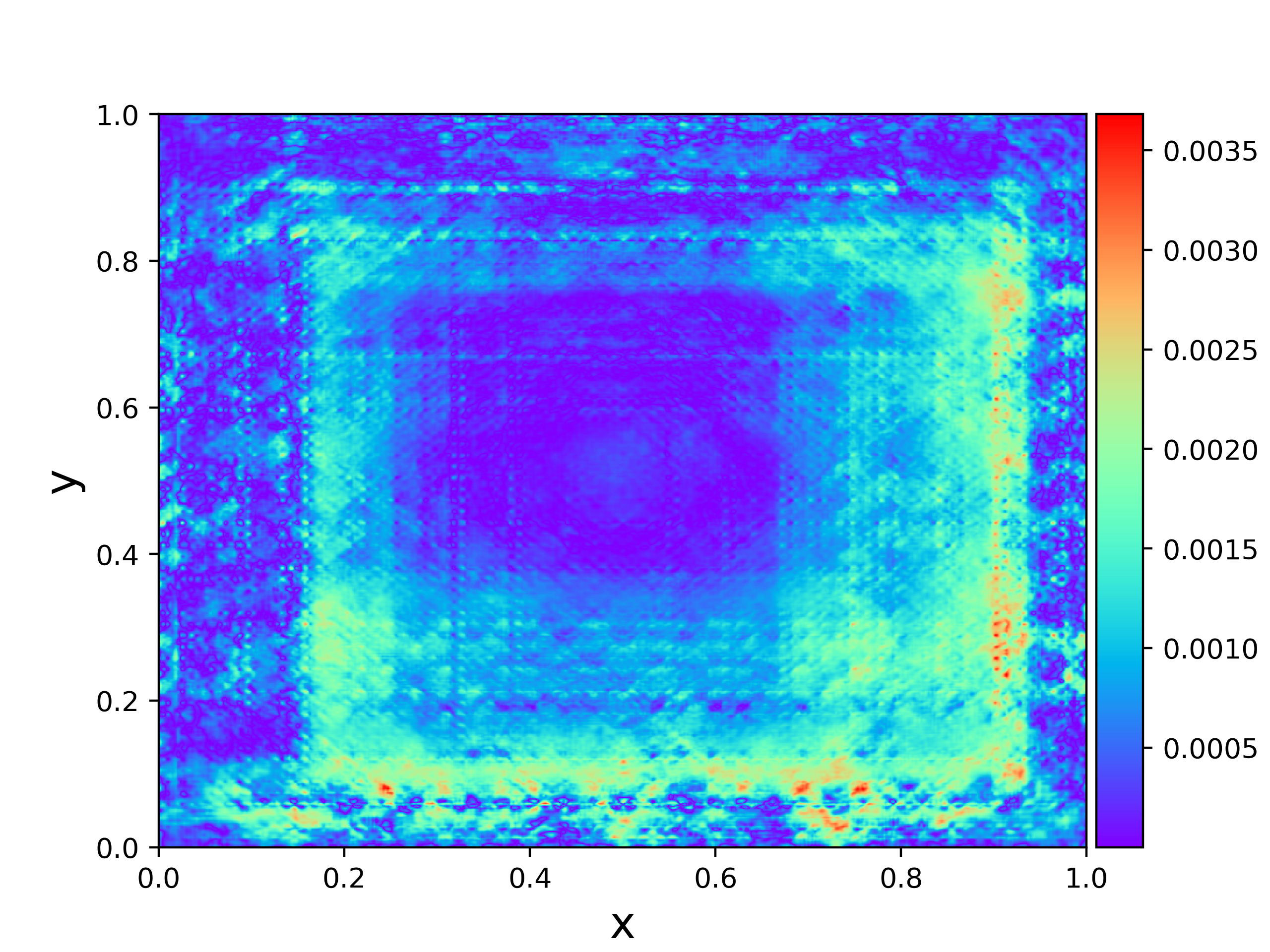}
	\caption{The absolute point-wise error for Example 7.}
	\label{figexam7}
\end{figure}

\begin{figure}[H]
	\centering
	\subfigure[Profile of $u(x,0.5)$]{ 
		\includegraphics[align=c,width=0.30\columnwidth]{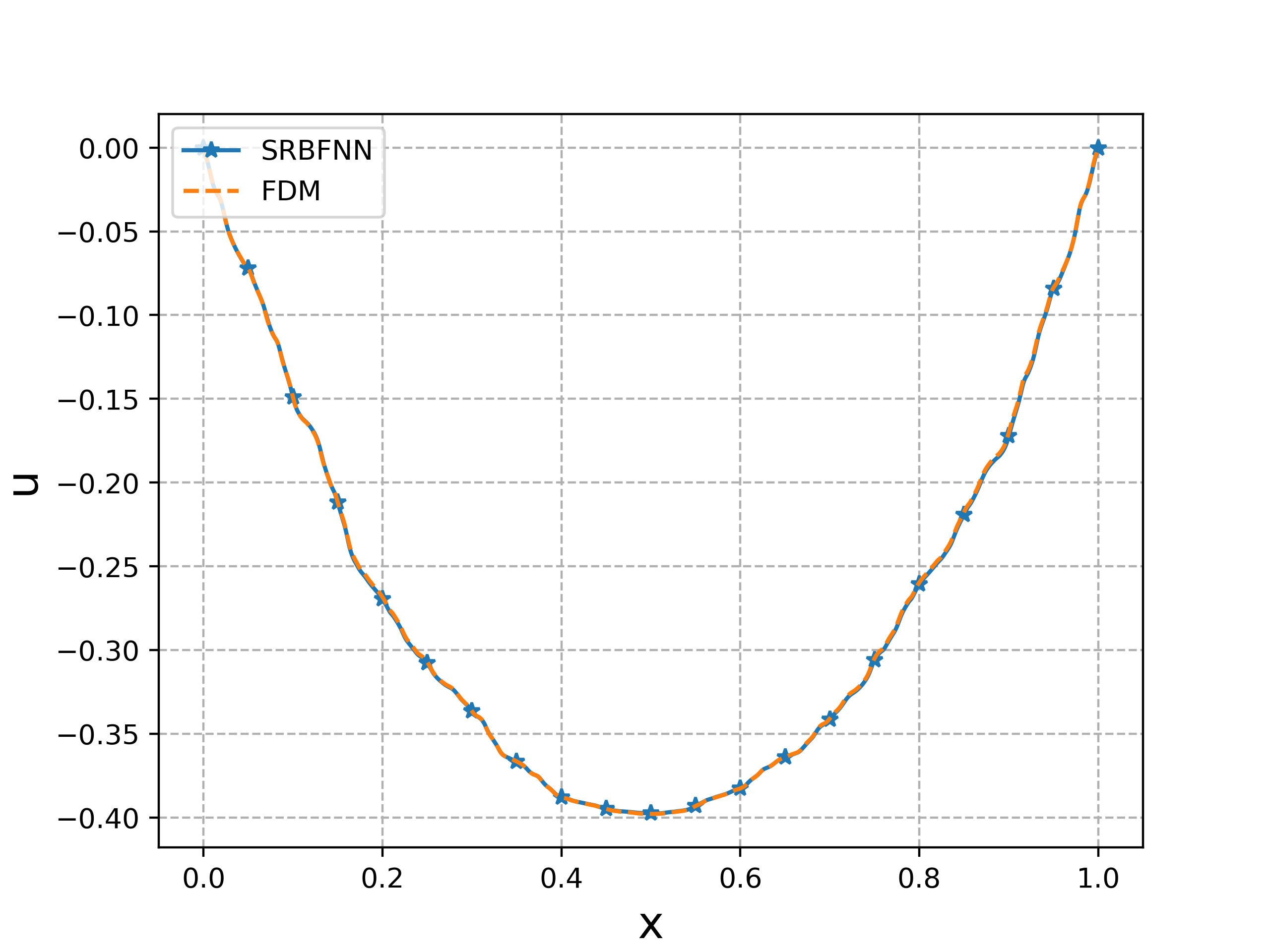}
	}
	\hfill %   ???   ?   
	\subfigure[Profile of $\frac{\partial u}{\partial x}|_{(x,0.5)}$]{ 
		\includegraphics[align=c,width=0.30\columnwidth]{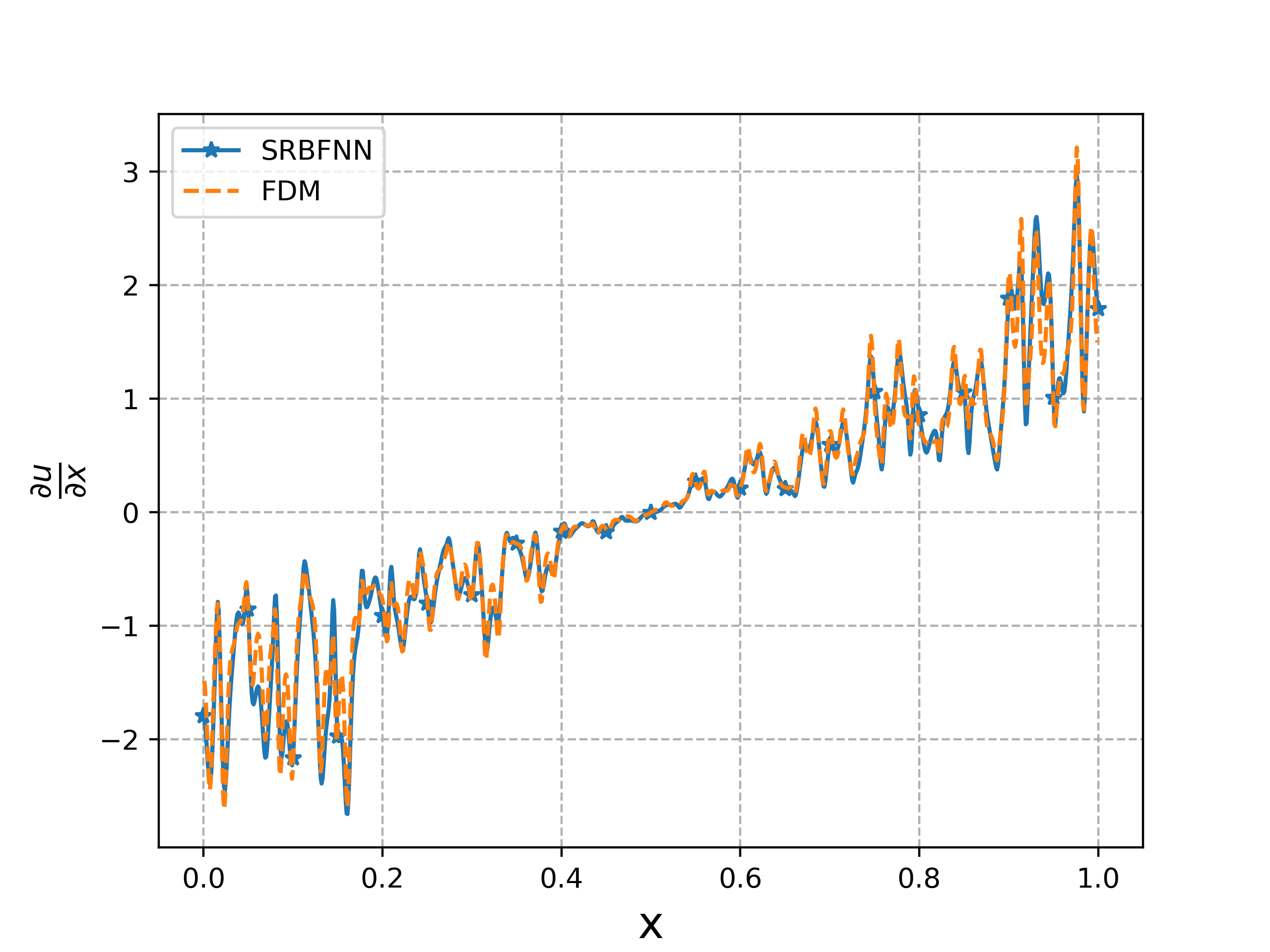}
	}
	\hfill %   ???   ?   
	\subfigure[Profile of $\frac{\partial u}{\partial y}|_{(0.5,y)}$]{ 
		\includegraphics[align=c,width=0.30\columnwidth]{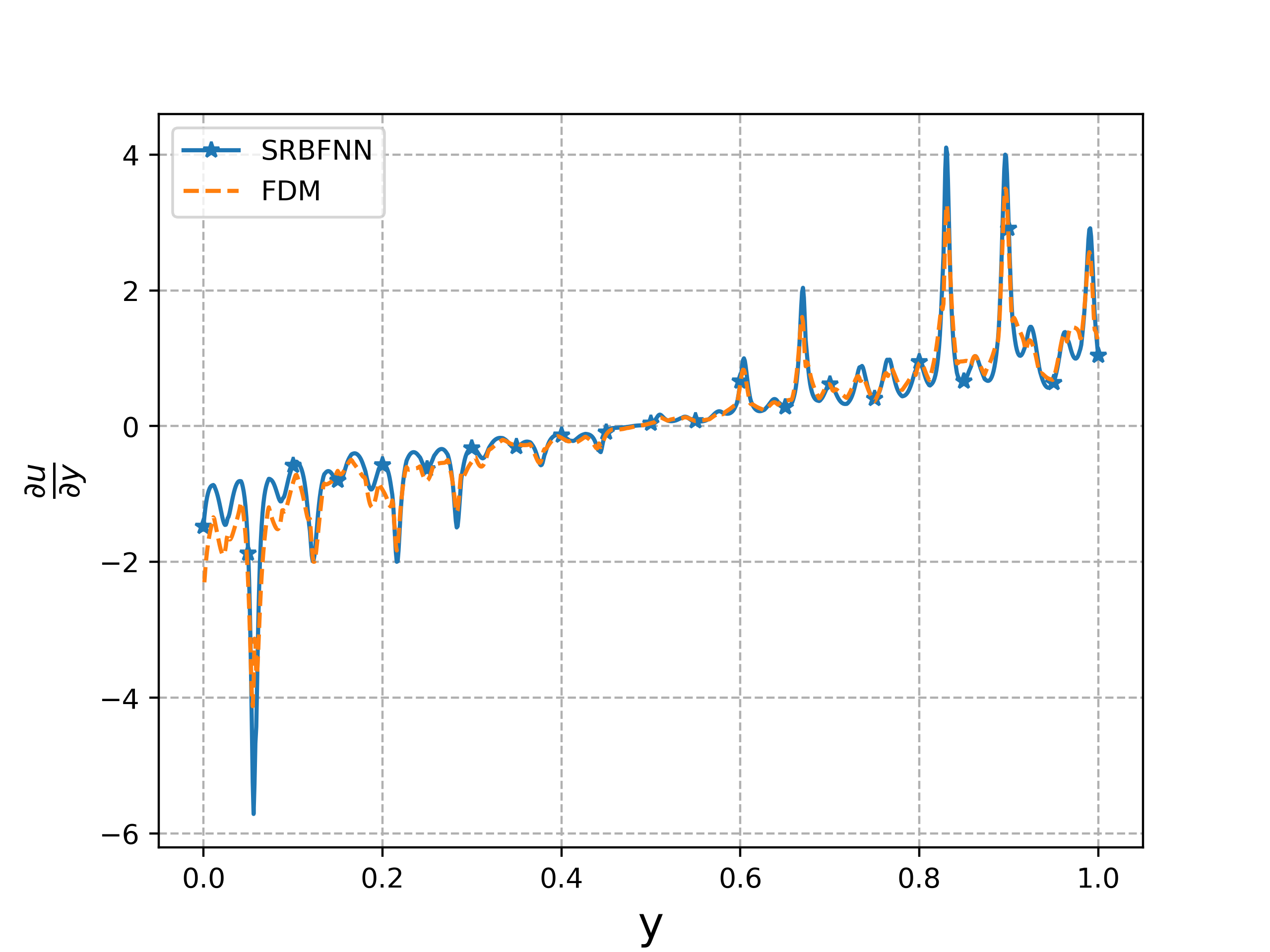}
	}
	\hfill %   ???   ?   
	\caption{The cross-section view of $u$, $\frac{\partial u}{\partial x}$ and  $\frac{\partial u}{\partial y}$ for Example 7. The dotted line represents the FDM solution. The solid line with an asterisk is the SRBFNN solution.}
	\label{fig11}
\end{figure}

The absolute point-wise error $|u^S-u^F|$ is shown in Fig. \ref{figexam7}.  SRBFNN can recover the solution up to $10^{-3}$ point-wise accuracy. To further analyze the performance of SRBFNN, we draw the cross-section view of the solution and its derivatives in Fig. \ref{fig11}. Table \ref{tabexam7} records three relative errors and the number of RBFs in the final solution. According to Fig. \ref{fig11} and Table \ref{tabexam7}, we can observe that SRBFNN achieves a good approximation in terms of the solution and its derivatives.
\begin{table}[H]
	\centering
	\begin{tabular}{cccc}
		\hline
		$N$   & $err_2$  & $err_\infty$ & $err_{H_1}$  \\ \hline	
		\rowcolor{gray!40}
		13213 & 4.427e-3 &1.249e-2      & 1.581e-1     \\ \hline
	\end{tabular}
	\caption{Relative errors and the number of basis functions in the final solution for Example 7.}
	\label{tabexam7}
\end{table}

%----------Numerical examples for three dimensions ----------------%
\subsection{Numerical experiment in 3D }
In this section, we consider a three-dimensional multiscale elliptic equation to verify the applicability of SRBFNN. In this scenario, solving multiscale elliptic equations with FDM is challenging. Thus we only consider the SRBFNN solution with a series of $\varepsilon=0.5, 0.2, 0.1, 0.05$ and take $u^{0.05}$ as the reference.

\begin{figure}[H]
	\centering
	\subfigure[$\varepsilon=0.5$]{ 
		\includegraphics[align=c,width=0.2\columnwidth]{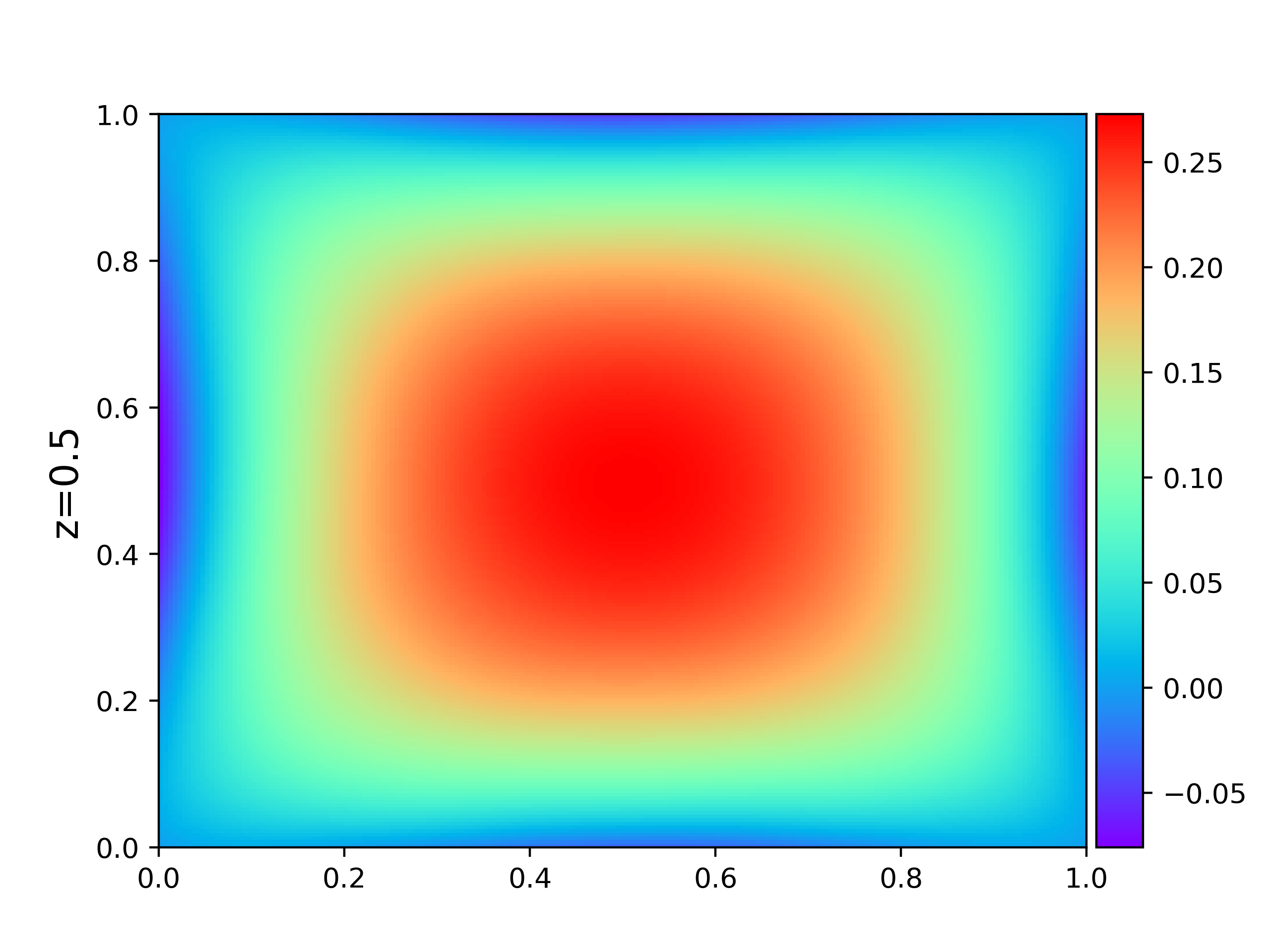}
	}
	\hfill %   ???   ?   
	\subfigure[$\varepsilon=0.2$]{ 
		\includegraphics[align=c,width=0.2\columnwidth]{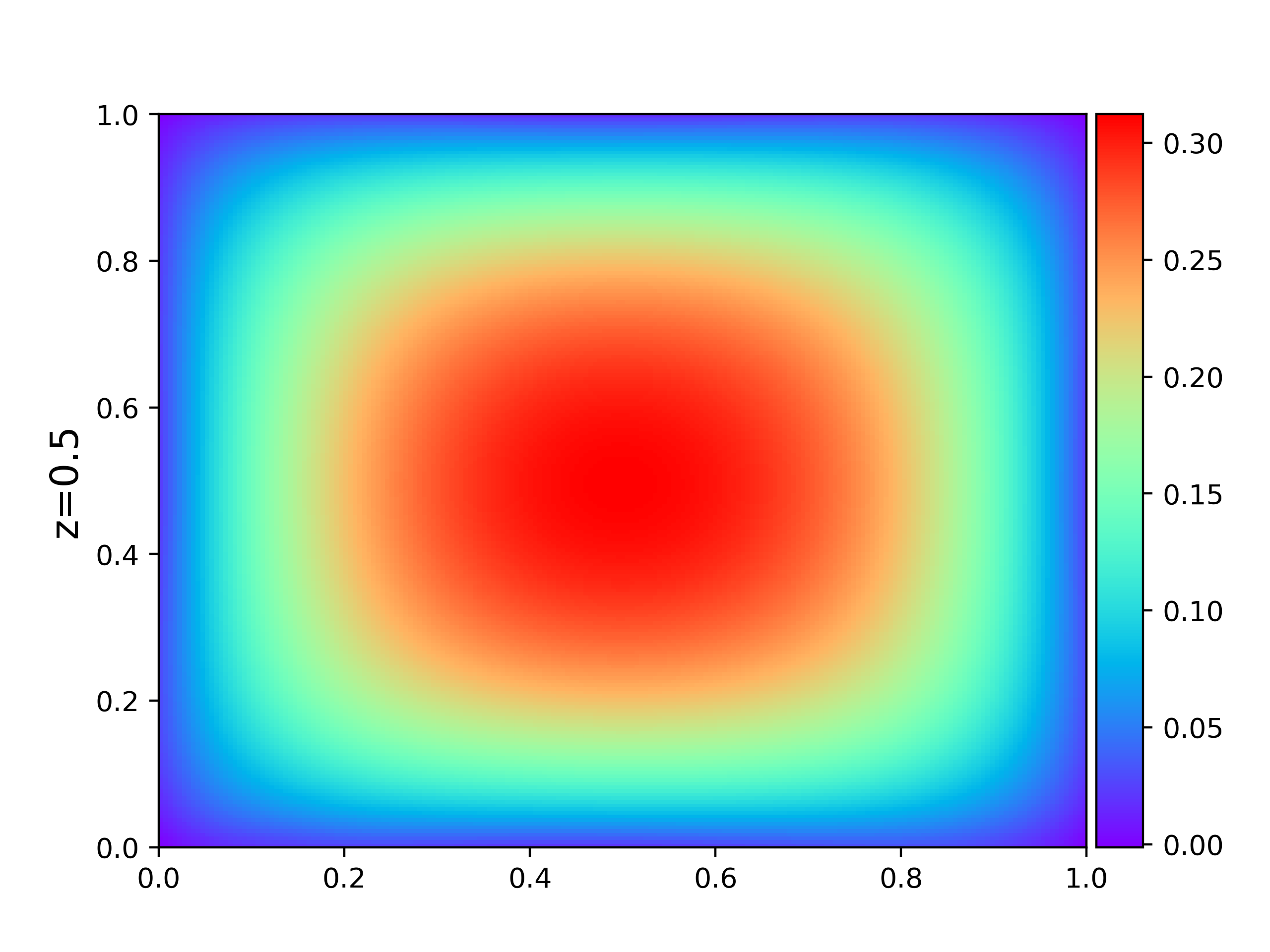}
	}
	\hfill %   ???   ?   
	\subfigure[$\varepsilon=0.1$]{ 
		\includegraphics[align=c,width=0.2\columnwidth]{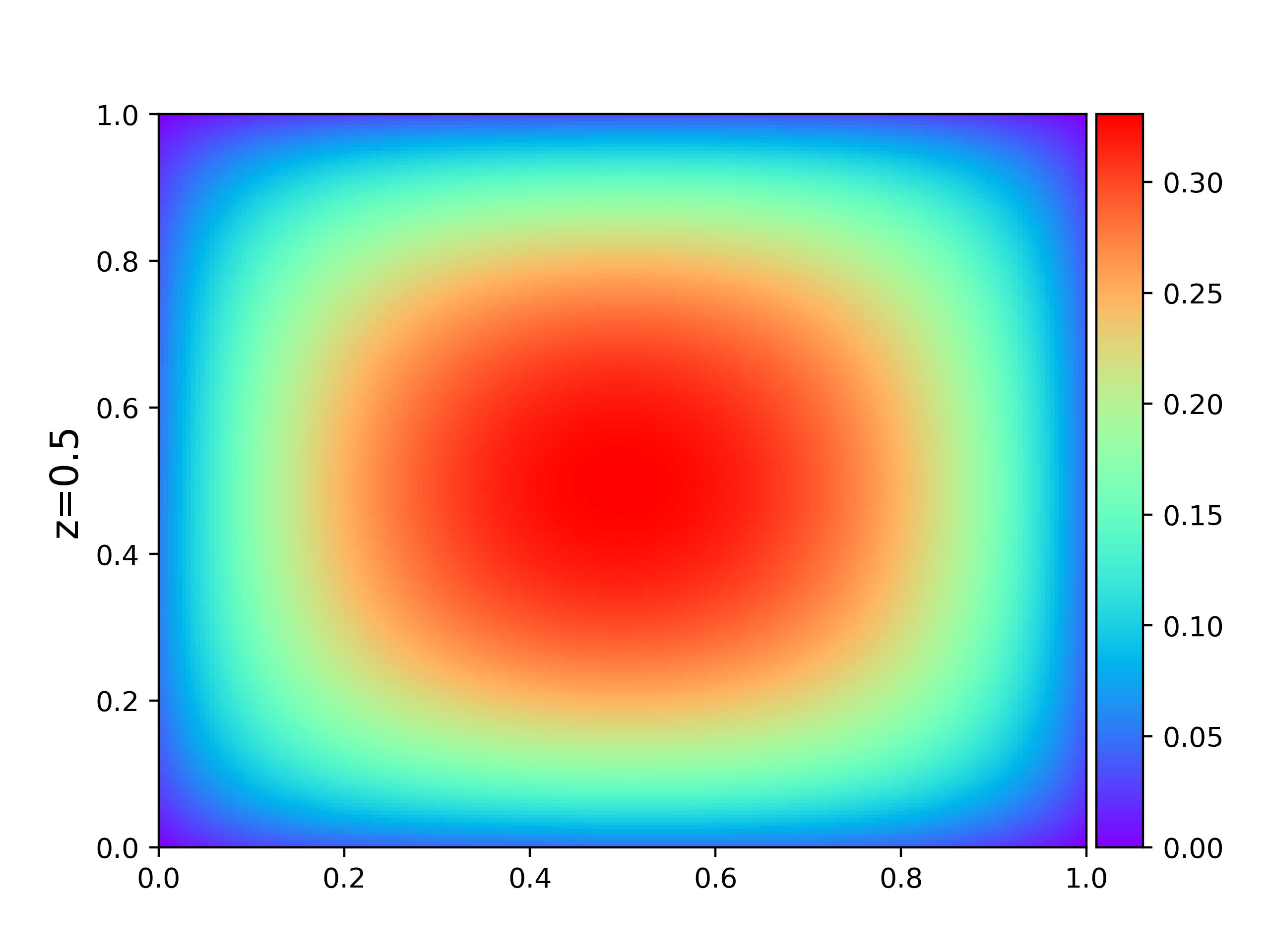}
	}
	\hfill %   ???   ?   
	\subfigure[$\varepsilon=0.05$]{ 
		\includegraphics[align=c,width=0.2\columnwidth]{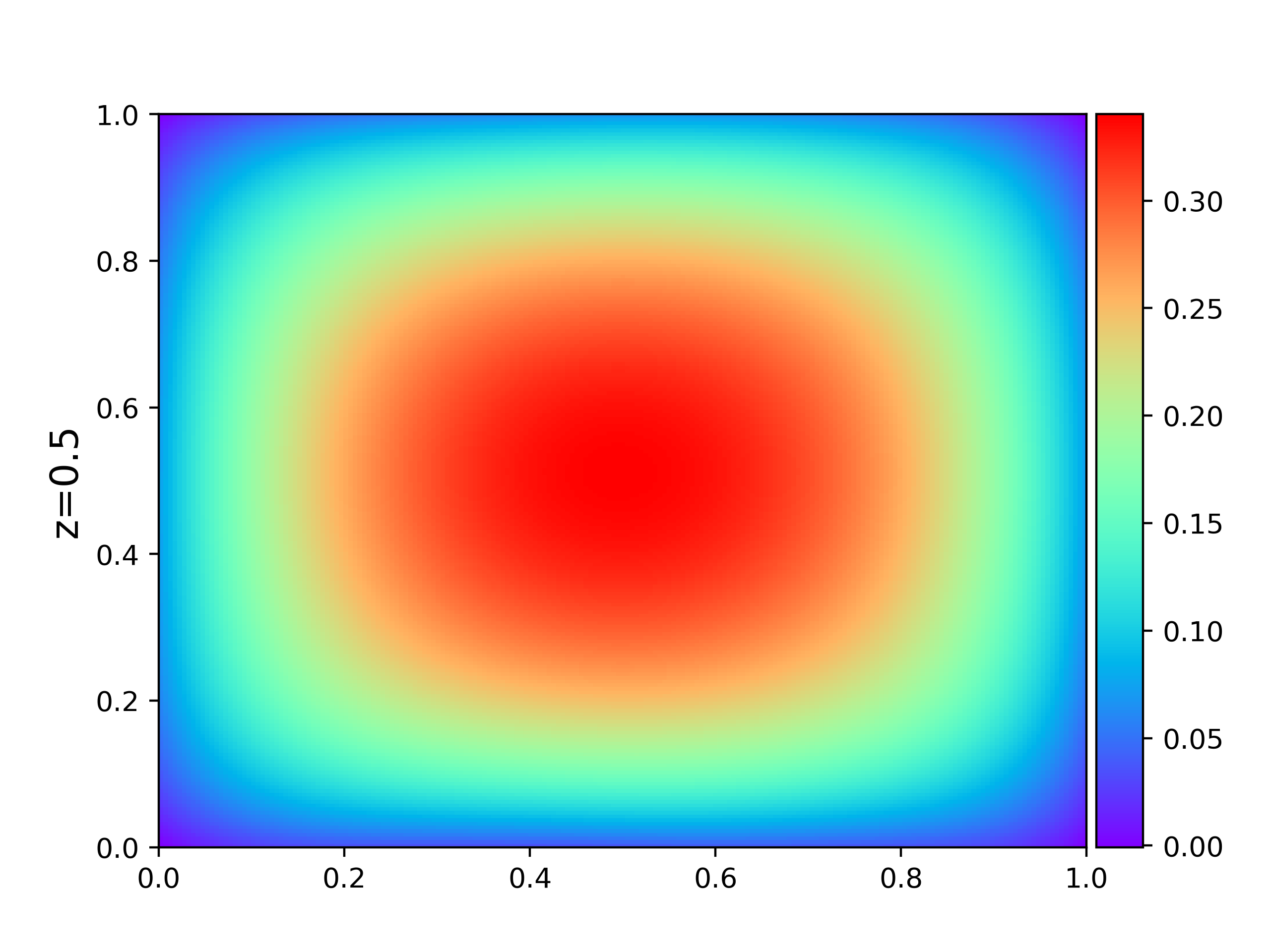}
	}
	\caption{Slices of the solution from SRBFNN with different $\varepsilon$ for Example 8 when $z=0.5$.}
	\label{figexam9}
\end{figure}

\begin{figure}[H]
	\centering
	\subfigure{ 
		\includegraphics[align=c,width=0.45\columnwidth]{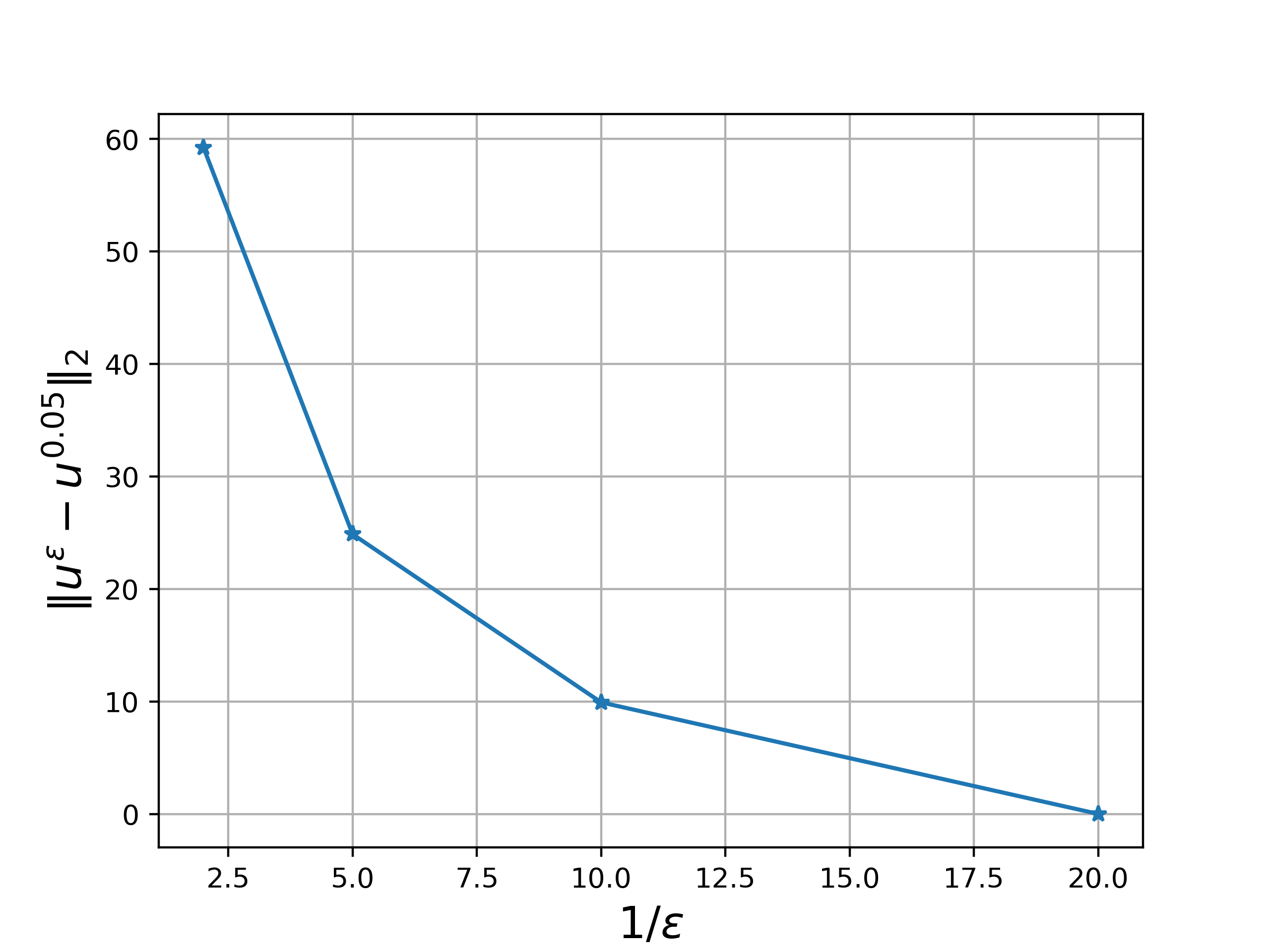}
	}
	\hfill %   ???   ?   
	\subfigure{ 
		\includegraphics[align=c,width=0.45\columnwidth]{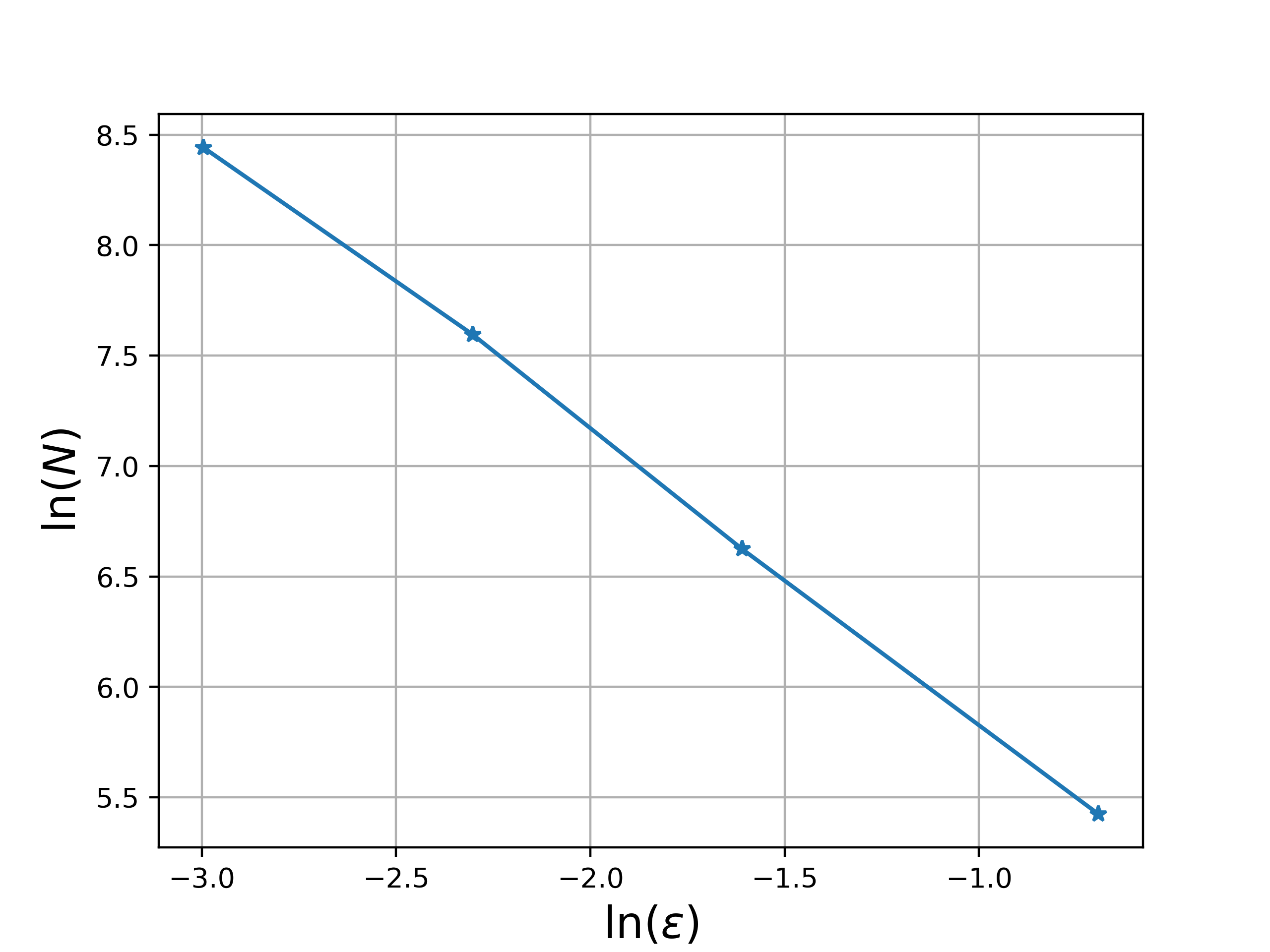}
	}
	\caption{The error between the solution when $\varepsilon=0.05$ and the solutions with larger scales and the relationship between the number of RBFs and $\varepsilon$. }
	\label{figexam9_1}
\end{figure}

\textbf{Setting:}\label{3dset} The training data is sampled equidistantly with mesh size $h$=0.01. $batchsize$ is 1024 in the domain and $200\times8$ on the boundary. $MaxNiter$ is 150 and $SparseNiter$ is 120. The initial learning rate is 0.1 and reduces by 1/10 every 30 iterations. The initial number of RBFs is 1000, 2000, 5000, 10000, respectively. $tol_1=0.05$, $tol_2=10^{-5}$, $Check_{iter}=10$. $\lambda_{1},\lambda_{2},\lambda_{3}$ are initialized to 0.1, 50.0, 0.001, respectively.
\paragraph{Example 8}
Consider the coefficient $a^{\varepsilon}(x,y,z)$ 
\begin{eqnarray}
	a^{\varepsilon}(x,y,z)=2+\sin(\frac{2\pi x}{\varepsilon})\sin(\frac{2\pi y}{\varepsilon})\sin(\frac{2\pi z}{\varepsilon}),
	\label{exam8}
\end{eqnarray}
where $f(x,y,z)=10.0$, $g(x,y)=0.0$. 
Fig. \ref{figexam9} shows the slices of the SRBFNN solution when $z=0.5$. We can see that the range of the solution is gradually the same as $\varepsilon$ tends to 0.05.
Table \ref{tabexam8} and Fig. \ref{figexam9_1} show the $L_2$ norm between $u^{0.05}$ and the solutions of other scales.
As expected, $u^{\varepsilon}$ approaches $u^{0.05}$ and the number of basis functions increase with the scale $\varepsilon$ tending to 0.05. These results  illustrate our method can achieve numerical convergence for three-dimensional multiscale elliptic equations.
\begin{table}[H]
	\centering
	\begin{tabular}{ccc}
		\hline
		$\varepsilon$ & $N$ & $\Vert u^\varepsilon-u^{0.05}\Vert_{L_2(\Omega)}$   \\ \hline	
		\rowcolor{gray!40}
		0.5	 &227  & 59.2  \\ 
		0.2	 &753  & 24.9  \\ 
		\rowcolor{gray!40}
		0.1  &1988 & 9.9  \\ 
		0.05 &4637 & 0.0  \\ \hline
	\end{tabular}
	\caption{The number of basis functions in the final solution and the $L_2$ error compared to the solution of SRBFNN when $\varepsilon$=0.05 in terms of $\varepsilon$ for Example 8.}
	\label{tabexam8}
\end{table}

%----------Conslusion----------------%
\section{Conslusion}\label{sec5}
This paper proposes a novel method to solve multiscale elliptic problems with oscillatory coefficients, such as the coefficients with scale separation, discontinuity and multiple scales. Unlike the general DNNs, SRBFNN has a simpler and more efficient architecture, which makes the network converge faster. The activation functions with good local approximation are important for multiscale problems. 
Furthermore, the $\ell_1$ regularization can deal with the overfitting problem owing to the simplicity of SRBFNN. The experiments show that SRBFNN presents the solution with fewer basis functions and has better approximation accuracy than most other deep learning methods. Finally, it is found that there is a relation between the number of basis functions $N$ in the final solution and the scale $\varepsilon$ in multiscale elliptic problems: $N= \mathcal{O}(\varepsilon^{-\tau n})$, $n$ is the dimensionality, and $\tau$ is typically smaller than $1$, which is better than classical numerical methods.

There are several interesting issues which deserve further consideration. First, the choice of the initial number of RBFs is essential, affecting the network's accuracy and training. Second, a wise selection of penalty parameters in different loss terms facilitates the training process and provides a better approximation. Another issue is to find a faster sparse optimization algorithm.

\section*{Acknowledgments}
This work is supported in part by National Key R$\&$D Program of China via grant No. 2022YFA1005203 and Jiangsu Provincial Key Research and Development Program under Grant BE2022058-4.

%%Vancouver style references.
\bibliographystyle{plain}
\bibliography{references.bib}

\end{document}